\documentclass{amsart}

\usepackage{stmaryrd}
\usepackage{xcolor}
\usepackage{tikz}
\usepackage{amsmath}
\usepackage{amsthm}
\usepackage{amssymb}
\usepackage{multicol}
\usepackage{hyperref}
\usepackage{graphicx}
\usepackage{caption} 
\usepackage{subcaption} 
\usepackage{geometry}
\usepackage{float}
\usepackage[ruled,vlined]{algorithm2e}

\numberwithin{equation}{section}
\newtheorem{theorem}{Theorem}[section]
\newtheorem{corollary}[theorem]{Corollary}
\newtheorem{example}[theorem]{Example}
\newtheorem{lemma}[theorem]{Lemma}
\newtheorem{proposition}[theorem]{Proposition}
\newtheorem{problem}{Problem}
\theoremstyle{definition}
\newtheorem{definition}[theorem]{Definition}
\newtheorem{remark}[theorem]{Remark}

\newcommand{\comments}[1]{}

\renewcommand{\textcolor}[2]{#2}
\title{Branching Pattern Modeling: A Geometric Measure Theory Approach}
\author{Qinglan Xia}
\address{
Department of Mathematics\\
University of California
at Davis\\
Davis, CA, 95616, USA }
\email{qlxia@ucdavis.edu}
\date{}


\usepackage{amsfonts}
\usepackage{amsmath,amssymb,amsthm}
\usepackage{graphicx, caption, subcaption}
\usepackage{tikz}
\usetikzlibrary{arrows.meta, positioning}
\usepackage{hyperref}

\numberwithin{equation}{section}

\theoremstyle{definition}

\newcommand{\Ccal}{\mathcal{C}}

\newcommand{\rect}{\underline{\underline{\tau}}}
\newcommand{\haus}[1][1]{\mathcal H^{#1}}
\newcommand{\chiT}[1][T]{\chi_{_{#1}}}
\newcommand{\malpha}[1][\alpha]{\mathbf M_{#1}}

\begin{document}

\begin{abstract}
\textcolor{blue}{This article is based on a lecture given at Westlake University. It aims to introduce students to an application of geometric measure theory in modeling branched patterns found in both living and nonliving systems. We begin with the Eulerian framework of ramified optimal transportation, including discrete transport paths, their basic properties, numerical approximation, and the passage to general transport paths represented by rectifiable currents. We then describe an application to the dynamic formation of a plant leaf, where branching venation patterns arise from a balance between metabolic benefit and transportation cost. Finally, we review landscape functions associated with transport paths and explain how two existing approaches can be unified within a broader family of rectifiable $1$-currents. We conclude by indicating several open directions related to Gilbert-type network problems, biological modeling, and numerical approximation.}
\end{abstract}
\date{}
\maketitle
\tableofcontents

\section{Introduction to ramified optimal transportation}
\subsection{Motivations from Plateau's problem.}
Nature presents an abundance of structures of remarkable beauty, many of which have inspired deep mathematical inquiry. One such example is the delicate geometry of soap films, whose forms motivated the study of the well-known \emph{Plateau's problem}: given a boundary curve, find a surface of least area spanning it. This question gave rise to the modern theory of minimal surfaces, and its solution depends critically on precise definitions of \emph{surface}, \emph{area}, and \emph{boundary}. Today, geometric measure theory (\cite{Federer1969, simon,Hardt2003size}) provides the rigorous framework for formulating and addressing such problems.

Another class of structures frequently encountered in nature comprises branching patterns. These occur across a wide spectrum of systems, both living and non-living, and likewise invite mathematical investigation. 
In biological contexts, branching structures are evident in trees, leaf venation, and the cardiovascular networks of animals. In non-biological systems, analogous patterns appear in river channel networks, railway and airline route systems, electric power distribution grids, and communication networks. The prevalence of these structures across diverse domains suggests that branching is often an efficient or optimal strategy for distribution and transport.

This ubiquity raises fundamental questions: Why are branching patterns favored by both natural processes and engineered systems? What advantages do they offer relative to non-branching configurations? Understanding the principles underlying these structures provides insight into efficiency, adaptability, and optimization in both natural and man-made networks.

\textcolor{blue}{In this article, we describe how geometric measure theory provides a natural framework for modeling
such branching structures. The central example is ramified optimal transportation, where transport paths
are represented by weighted graphs or, more generally, by rectifiable $1$-currents, and where branching arises
from economies of scale in the transportation cost. This viewpoint connects network optimization, variational
problems, and geometric measure theory.}

\textcolor{blue}{For broader background on branched transport, traffic plans, irrigation patterns, and related variational
formulations, we refer the reader to \cite{MSM2003, BCM2005, bernot2008optimal, BBS2006, BBS2011, MS2010}.
For related geometric-measure-theoretic and variational perspectives on branched transport, including
flat-chain formulations, size minimization, stability, and well-posedness, see
\cite{Hardt2003size, paolini2006optimal, bw, colombo2018, colombo2020}.}

\textcolor{blue}{We now begin with the simplest setting: ramified optimal transportation between atomic measures (\cite{xia2003, xia2015motivations}).}

\subsection{Discrete version of transport paths}
Let us begin with a simple motivating example. Suppose two products are to be transported from nearby cities 
to a distant city. One option is to ship them separately, resulting in a “V-shaped” path; another is to first bring them to a common location and then transport them together in a single shipment, producing a “Y-shaped” path. In many situations, the latter is preferred due to \emph{economies of scale} in transportation, as combining shipments for a long trip often reduces overall costs.

This example illustrates a fundamental principle in network design: merging flows at intermediate junctions can minimize overall transport effort, a concept that directly motivates the study of \emph{transport paths} in both natural and engineered systems.

\subsubsection{Set-up}
Let $X$ be a convex compact subset of a Euclidean space $\mathbb{R}^{m}$. For any $x \in X$, let $\delta_x$ denote the Dirac measure centered at $x$. A finitely supported atomic measure in $X$ is of the form
\[
\sum_{i=1}^{k} m_i \delta_{x_i},
\]
with distinct points $x_i \in X$ and $m_i > 0$ for each $i = 1, \dots, k$. The value $\Lambda = \sum_{i=1}^k m_i$ is called the total mass of the measure. For any $\Lambda > 0$, let
\begin{equation}
\mathcal{A}_{\Lambda }\left( X\right)
\label{a_lambda}
\end{equation}%
be the space of all atomic measures on $X$ with total mass $\Lambda$.

\begin{definition}
    
For any $\Lambda >0$, and any two atomic measures%
\begin{equation}
\mathbf{a}=\sum\limits_{i=1}^{k}m_{i}\delta _{x_{i}}\text{ and }\mathbf{b}%
=\sum\limits_{j=1}^{\ell }n_{j}\delta _{y_{j}}  \label{atomicmeasure}
\end{equation}%
in $\mathcal{A}_{\Lambda }\left( X\right) $, a \textit{transport path} from $%
\mathbf{a}$ to $\mathbf{b}$ is a weighted directed graph $G = \left(V(G), E(G), w\right)$, where
\begin{itemize}
\item $V(G) \subseteq X$ is the vertex set, representing the geometric location of vertices,
\item $E(G)$ is the set of directed edges, and
\item $w : E(G) \to (0, +\infty)$ assigns a positive weight to each edge
\end{itemize}
such that 
the vertex set contains all source and target points:
\[
\{x_1, \dots, x_k\} \cup \{y_1, \dots, y_\ell\} \subset V(G),
\]
and for each vertex $v \in V(G)$, the following \emph{balance equation} holds:
\begin{equation}
\sum_{\substack{ e\in E\left( G\right)  \\ e^{-}=v}}w\left( e\right) =\sum 
_{\substack{ e\in E\left( G\right)  \\ e^{+}=v}}w\left( e\right) +
\left\{ 
\begin{array}{c}
m_{i},\text{\ if }v=x_{i}\text{\ for some }i=1,...,k \\ 
-n_{j},\text{\ if }v=y_{j}\text{\ for some }j=1,...,\ell \\ 
0,\text{\ otherwise. }%
\end{array}%
\right.
  \label{balance}
\end{equation}
Here, $e^-$ and $e^+$ denote the starting and ending vertices of the edge $e \in E(G)$, respectively.  
\end{definition}

The balance equation above represents the conservation of mass at each vertex. In other words, $G$ satisfies
Kirchhoff's law at each vertex. 

In terms of polyhedral chains, a transport path $G$ from $%
\mathbf{a}$ to $\mathbf{b}$ is simply a real polyhedral 1-chain 
\begin{equation*}
G=\sum_{e\in E\left( G\right) }w\left( e\right) \left[ e\right]
\end{equation*} with $\partial G = \mathbf{b} - \mathbf{a}$.

\begin{remark}

\label{remark:path}
Each transport path $G$ corresponds to a vector-valued
measure on $X$
\begin{equation*}
G=\sum_{e\in E\left( G\right) }w\left( e\right) \left[ \left[ e\right] %
\right] ,
\end{equation*}
where $\left[ \left[ e\right] \right] $ is the vector-valued measure $\mathcal{H}%
^{1}\lfloor _{e}\overrightarrow{e}$ for each edge $e\in E\left( G\right) $
with unit directional vector $\overrightarrow{e}$. The boundary conditions can
be simplified to be a single divergence condition:
\begin{equation*}
\text{div}\left( G\right) =\mathbf{a-b},
\end{equation*}
in the sense of distribution. 
\end{remark}

Finally, let
\[
Path(\mathbf{a,b})
\]
denote the space of all transport paths from $\mathbf{a}$ to $\mathbf{b}$. Among all paths in $Path(\mathbf{a,b}) $, we want to find an optimal path
which allows the possibility that some parts overlap in a cost efficient
fashion. To get such a \textquotedblleft Y-shaped\textquotedblright\ optimal
path, we define the following cost function on transport paths.

\begin{definition}

For any $\alpha \in \left[ 0,1\right] $, the $\mathbf{M}_{\alpha }$\textbf{\ 
}cost of a transport path \[G=\left\{ V\left( G\right) ,E\left( G\right)
,w:E\left( G\right) \rightarrow \left( 0,+\infty \right) \right\} \] is
defined by 
\begin{equation}
\mathbf{M}_{\alpha }\left( G\right) :=\sum_{e\in E\left( G\right) }w\left(
e\right) ^{\alpha }length\left( e\right) ,  \label{alpha_cost}
\end{equation}%
where $length(e)$ denotes the length of the edge $e$, typically the Euclidean distance between endpoints $e^-$ and $e^+$ of $e$.

\end{definition}

We now consider the following  problem:

\begin{problem}[Ramified optimal transport problem: Discrete Version] Given two finitely supported atomic measures $\mathbf{a,b}\in \mathcal{A}_{\Lambda }\left( X\right)$ and $\alpha\in[0,1)$, minimize $\mathbf{M}_{\alpha}(G)$ among all transport paths $G$ in $Path(\mathbf{a,b})$.
\label{problem 1}
\end{problem}
\noindent In other words, it is a Plateau-type problem: minimize $\mathbf{M}_{\alpha}(G)$ among all real polyhedral 1-chains $G$ with $\partial G=\mathbf{b-a}$.

An $%
\mathbf{M}_{\alpha }$ minimizer in $Path(\mathbf{a},\mathbf{b})$ is called
an $\alpha -$\textit{optimal transport path} from $\mathbf{a}$ to $\mathbf{b}
$.

\begin{remark}
    The parameter $\alpha\in[0,1]$ is used to measure the transport economy of scale: When $\alpha<1$, the function $x^\alpha$ is concave, the inequality $(m_1+m_2)^\alpha<m_1^\alpha+m_2^\alpha$ makes the overall cost cheaper to transport products in groups rather than separately. Also, when $\alpha=1$, it corresponds to Monge optimal transport problem\cite{villani2009book}, while $\alpha=0$ is analogous to the Steiner tree problem. Since the case $\alpha=1$ does not exhibit branching phenomena, we usually only consider $\alpha<1$ in the ramified optimal transportation.
\end{remark}

\begin{remark}
\textcolor{blue}{Gilbert's minimum-cost communication network problem \cite{Gilbert1967} is a classical precursor to ramified optimal transportation. It already captures the key economy-of-scale principle: carrying many channels along a shared route may reduce the cost per unit. The formulation used here, however, is the modern transport-path formulation with directed weighted graphs, prescribed boundary $\partial G=b-a$, and cost $\mathbf M_\alpha$. In particular, it may be viewed as a Plateau-type problem for weighted polyhedral $1$-chains; as discussed below, optimal transport paths can be taken to be acyclic. Gilbert's original problem was formulated instead as a communication-network design problem with pairwise channel demands between cities, where the network links are undirected or bidirectional rather than directed transport edges. Closely related formulations are often discussed under the name ``Gilbert problem'' in the later branched-transport literature, especially in \cite{bernot2008optimal}.}
\end{remark}

\begin{remark}
When the author introduced transport paths in \cite{xia2003} using an Eulerian approach, a Lagrangian model was introduced essentially at the same time by Maddalena, Solimini and Morel \cite{MSM2003}, and then intensively studied by Bernot, Caselles, Morel \cite{BCM2005, bernot2008optimal}. In this approach, each transport path is viewed as a measure on the space of Lipschitz curves. We refer to the book \cite{bernot2008optimal} for the equivalence of these approaches. In Section \ref{sec: Santambrogio}, we will recall some basic definitions of this approach and introduce landscape functions there.
\end{remark}

\subsubsection{Basic properties of optimal transport paths}
We now describe some basic properties of transport paths as studied in \cite{xia2003}.

Note that an arbitrary transport path $G\in Path\left( \mathbf{a,b}\right) $ is a weighted
directed graph, but not necessarily a directed tree. In other words, the support of $G$ may
contain some (undirected) cycles. Here, a weighted directed graph $G=\left\{ V\left(
G\right) ,E\left( G\right) ,w:E\left( G\right) \rightarrow (0,1]\right\} $
contains a \textit{cycle} if for some $k\geq 3$, there exists a list of
distinct vertices $\left\{ v_{1},v_{2},\cdots ,v_{k}\right\} $ in $V\left(
G\right) $ such that for each $i=1,\cdots ,k$, either\ the segment $\left[
v_{i},v_{i+1}\right] $ or $\left[ v_{i+1},v_{i}\right] $ is a directed edge
in $E(G)$, with the agreement that $v_{k+1}=v_{1}$. However, the following
proposition says $G$ can be modified to be an acyclic graph $\tilde{G}\in
Path\left( a,b\right) $  with less or equal $\mathbf{M}%
_{\alpha }$ cost.

\begin{proposition}
\cite[Proposition 2.1]{xia2003}
\label{loop} Given $\mathbf{a,b}$ as in (\ref{atomicmeasure}) and $\alpha\in [0,1)$. For any transport path $G\in Path\left( \mathbf{a,b}\right) $,
there exists a transport path $\tilde{G}\in Path\left( \mathbf{a,b}\right) $
such that $V\left( \tilde{G}\right) \subseteq V\left( G\right) $, $\mathbf{M}%
_{\alpha }\left( \tilde{G}\right) \leq \mathbf{M}_{\alpha }\left( G\right) $
and the support of $\tilde{G}$ contains no (undirected) cycles.
\end{proposition}
Let 
\begin{equation*}
Path_{0}\left(\mathbf{ a,b}\right) =\left\{ G\in Path\left( \mathbf{a,b}\right): \text{ the support of } G \text{
contains no cycles}\right\} .
\end{equation*}
When considering Problem \ref{problem 1}, from the above proposition, we may restrict our transport paths to elements in $Path_{0}\left(\mathbf{ a,b}\right)$. The following proposition says that the number of all branching vertices of $G\in Path_{0}\left(\mathbf{ a,b}\right) $ is bounded above.
\begin{proposition}
\label{prop: upperbound}
Suppose $G\in Path_{0}\left(\mathbf{ a,b}\right) $. Then, 
\begin{equation*}
\left\vert \left\{ v:\deg \left( v\right) \geq 3\right\} \right\vert \leq
k+\ell -2,
\end{equation*}%
where $deg(v)$, called the degree at $v$, is the number of edges in $E(G)$ having $v$ as one of its endpoints. 
\end{proposition}
\begin{proof}
By calculating the total number of vertices in $V(G)$, we have
\[2|E(G)|=\sum_{v\in V(G)} deg(v).\]
By means of the Euler number $\chi_G=|V(G)|-|E(G)|$ of $G$, we have 
\[2\sum_{v\in V(G)}1-2\chi_G=2|V(G)|-2\chi_G=2|E(G)|=\sum_{v\in V(G)} deg(v).\]
That is,
\[\sum_{v\in V(G), deg(v)=1}1-2\chi_G=\sum_{v\in V(G), deg(v)\ge 3} (deg(v)-2).\]
Hence,
\[\left\vert \left\{ v:\deg \left( v\right) \geq 3\right\} \right\vert  \le \sum_{deg(v)\ge 3} (deg(v)-2)=\sum_{ deg(v)=1}1-2\chi_G\le k+\ell -2.\]
\end{proof}
As a result, given $k$ and $\ell $, there are finitely many possibilities of
topologies of all transport paths $G$ in $Path_{0}\left( \mathbf{a,b}\right) 
$. For instance, when $k=2$, $\ell =1$, \textcolor{blue}{there are four possible configurations: Y-shaped, V-shaped, 7-shaped, and reversed 7-shaped.} It is clear that there exists at least one
minimizer in each given topology. As a result, there always exists an $%
\mathbf{M}_{\alpha }$-minimizer in $Path\left( \mathbf{a,b}\right) $. That is, there
exists an $\alpha -$optimal transport path in $Path_{0}\left( \mathbf{a,b}%
\right) $ for any $\alpha\in [0,1)$.

\subsection{Numerical approximation of optimal transport paths}
\begin{figure}[htbp]
    \centering
    \includegraphics[width=\textwidth]{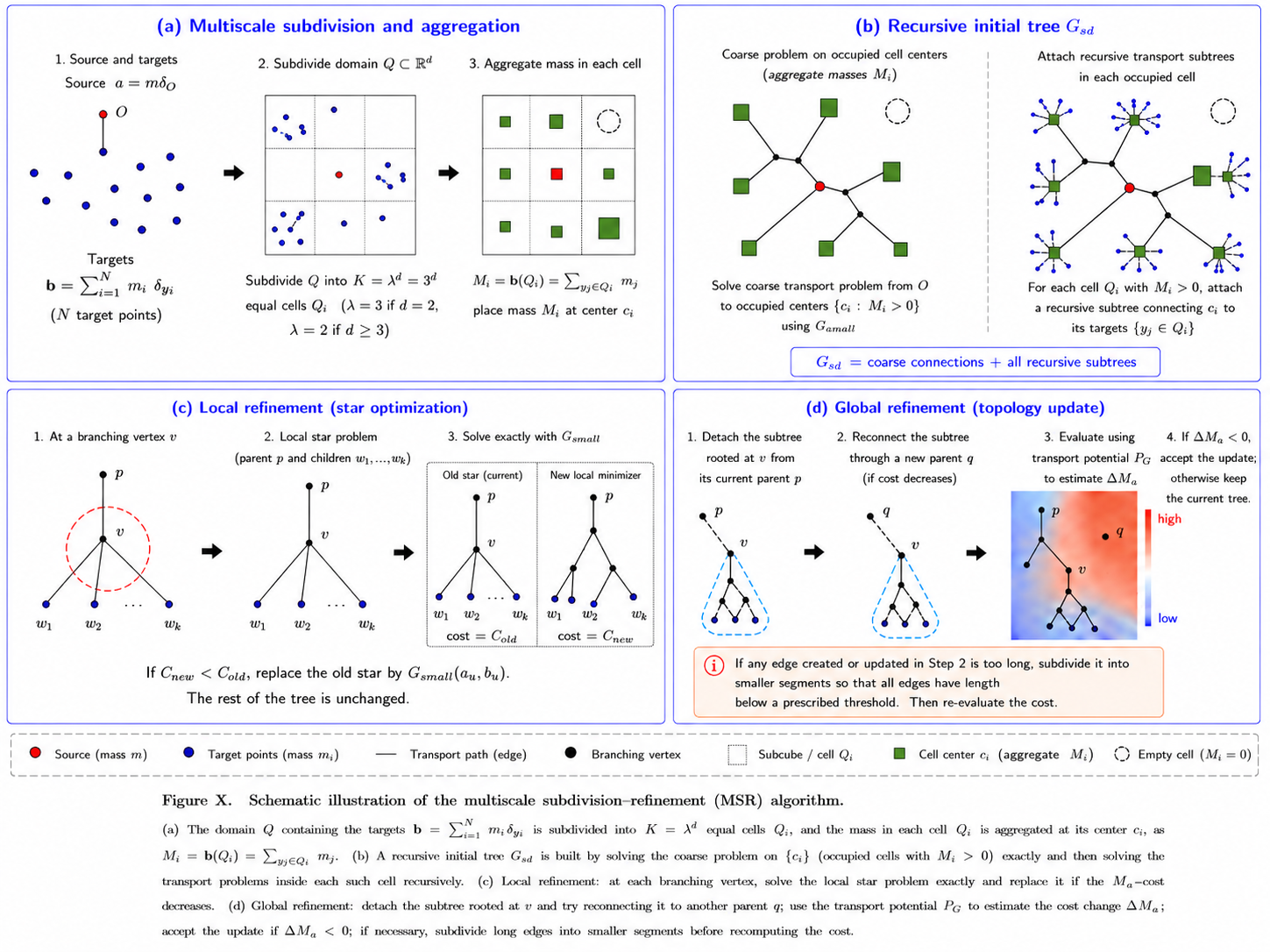}
    \caption{Schematic illustration of the multiscale subdivision--refinement (MSR) algorithm. 
Panel (a) shows the subdivision of the domain and the aggregation of target masses at cell centers. 
Panel (b) shows the recursive construction of the initial subdivision tree \(G_{sd}\). 
Panel (c) illustrates the local refinement step by star optimization. 
Panel (d) illustrates the global refinement step, in which a subtree is reassigned to a new parent, with long edges subdivided if necessary.}
\label{fig:MSR-algorithm}
\end{figure}

While optimal transport paths can be computed exactly for small atomic measures, numerical methods are needed to treat larger problems. Since the number of admissible transport tree topologies grows rapidly with the number of target points, exhaustive optimization soon becomes impractical. To visualize optimal transport paths numerically, we employ a multiscale subdivision--refinement (MSR) algorithm developed in \cite{xia2015motivations}, which combines exact optimization on small subproblems with iterative local and global refinements. The main steps of the algorithm are summarized schematically in Figure~\ref{fig:MSR-algorithm}.

The algorithm begins with an exact solver for small transport problems. Given atomic measures
\[
\mathbf a = m\delta_O,
\qquad
\mathbf b = \sum_{i=1}^{N} m_i\delta_{y_i},
\]
when the number of terminals is sufficiently small, all admissible transport tree topologies can be enumerated. For each fixed topology, the branching locations are optimized using a modified Melzak algorithm \cite{melzak}. Among all topologies, the one with the smallest \(M_\alpha\)-cost is then selected. We denote the resulting exact minimizer by
\[
G_{\mathrm{small}}(\mathbf a,\mathbf b).
\]

For larger problems, the main idea is to first build a reasonable transport tree by a coarse-to-fine spatial subdivision. Let
\[
\lambda =
\begin{cases}
3, & d=2,\\
2, & d\geq 3,
\end{cases}
\qquad
K=\lambda^d,
\]
where \(d\) is the ambient dimension. If \(N\leq K\), we set
\[
G_{sd}(\mathbf a,\mathbf b)
=
G_{\mathrm{small}}(\mathbf a,\mathbf b).
\]
If \(N>K\), let \(Q\) be a cube containing the support of \(\mathbf b\), and subdivide \(Q\) into \(K\) equal subcubes \(Q_1,\ldots,Q_K\). For each subcube \(Q_i\), let
\[
M_i
=
\mathbf b(Q_i)
=
\sum_{y_j\in Q_i} m_j
\]
be the total mass contained in \(Q_i\), and let \(c_i\) be the center of \(Q_i\). The masses inside \(Q_i\) are replaced at the coarse level by the aggregated mass \(M_i\delta_{c_i}\). The initial subdivision tree is then defined recursively by
\[
G_{sd}(\mathbf a,\mathbf b)
=
G_{\mathrm{small}}
\left(
\mathbf a,
\sum_{i=1}^{K} M_i\delta_{c_i}
\right)
+
\sum_{i=1}^{K}
G_{sd}
\left(
M_i\delta_{c_i},
\mathbf b\!\restriction_{Q_i}
\right).
\]
Thus \(G_{sd}\) first connects the source to a finite number of cluster centers and then recursively resolves the transport structure inside each cluster.

Starting from \(G_{sd}\), the MSR algorithm improves the transport tree through two refinement procedures, illustrated in panels (c) and (d) of Figure~\ref{fig:MSR-algorithm}. The first is a local optimization step. Around each branching vertex, the incident edges determine a small star-shaped transport problem involving the parent of the vertex and its children. Whenever the exact solver \(G_{\mathrm{small}}\) finds a lower-cost replacement for this local configuration, the corresponding subnetwork is replaced. This step decreases the transport cost locally while preserving the rest of the tree.

The second refinement is a global topology update. Local optimization alone may not change the large-scale branching structure. To overcome this limitation, the algorithm allows a subtree rooted at a vertex \(v\) to be detached from its current parent and reconnected through a different parent whenever doing so decreases the transport energy. This reassignment is guided by a transport potential \(P_G\), which estimates the cost variation caused by modifying the flow along the tree. More precisely, for a vertex \(u\) and a mass variation \(t\), one defines \(P_G(u,t)\) recursively by
\[
P_G(O,t)=0,
\]
and, for \(u\neq O\),
\[
P_G(u,t)
=
P_G(p(u),t)
+
|p(u)-u|
\left[
m(u)^\alpha
-
\bigl(m(u)-t\bigr)^\alpha
\right],
\]
where \(p(u)\) is the parent of \(u\), and \(m(u)\) is the mass carried by the edge from \(p(u)\) to \(u\). This potential provides an efficient way to compare possible parent reassignments. If necessary, long edges created during the global refinement step are subdivided into smaller segments before the cost is re-evaluated. In this way, the algorithm can perform nonlocal topological corrections that are not accessible through purely local star replacements.

\begin{remark}
The refinement procedures described above are not limited to the subdivision tree \(G_{sd}\). 
Whenever an initial transport path is already available, one may apply the local and global refinement steps directly to this path in order to further decrease its \(M_\alpha\)-cost.
\end{remark}

The computational cost of the MSR algorithm comes from the interaction between exact small-scale optimization and large-scale refinement. The exact solver $G_{\mathrm{small}}$ requires enumeration of admissible tree topologies, whose number grows combinatorially with the number of terminals. For this reason, it is used only on subproblems of bounded size. In the subdivision stage, $K=\lambda^d$ is fixed once the ambient dimension $d$ is fixed, so each call to $G_{\mathrm{small}}$ involves only a bounded number of aggregated target masses. The full problem is therefore replaced by a hierarchy of many small exact problems. For reasonably balanced point distributions, the depth of the subdivision tree grows logarithmically with the number $N$ of target points, and the construction of the initial subdivision tree $G_{sd}$ has practical complexity close to $O(N\log N)$, often behaving nearly linearly in numerical experiments.

After the initial tree has been constructed, the refinement stage usually dominates the running time. The local refinement step examines star-shaped neighborhoods around branching vertices. For an $\alpha$-optimal transport path with $N$ target points, Proposition~\ref{prop: upperbound} implies that the number of vertices is bounded above by $2N$. Moreover, the degree of each vertex in an optimal transport path is uniformly bounded above by a constant depending only on the model parameters. These structural properties explain why local refinement is computationally feasible: once the tree is close to an optimal configuration, each local exact problem involves only a bounded number of terminals, and one local refinement sweep has cost of order $O(N)$. For a general initial transport path, high-degree vertices may occur before refinement. In that case, the cost of one local sweep can be written schematically as
\[
\sum_{v\in V(G)} C_{\mathrm{small}}(\deg(v)+1),
\]
where $C_{\mathrm{small}}(k)$ denotes the cost of solving an exact $k$-terminal subproblem.

The global refinement step is more expensive, since it considers whether a subtree rooted at a vertex should be detached from its current parent and reconnected through another parent. Since a tree-like optimal transport path has at most $2N$ vertices, a naive implementation that compares all possible parent reassignments has worst-case cost of order $O(N^2)$ per global sweep. In practice, this cost can be reduced by restricting the search to geometrically relevant candidate parents and by using spatial search structures. The subdivision of long edges may introduce auxiliary degree-two vertices and thereby increase the number of numerical vertices, but these additional vertices are inserted for geometric resolution and can be controlled separately from the combinatorial branching structure.

Combining these estimates, if $r_{\mathrm{loc}}$ and $r_{\mathrm{glob}}$ denote the numbers of local and global refinement sweeps, respectively, then a naive implementation has total complexity of order
\[
O\bigl(N\log N+r_{\mathrm{loc}}N+r_{\mathrm{glob}}N^2\bigr).
\]
When the number of refinement sweeps is moderate, the global reassignment step is the dominant term in this estimate. With geometrically restricted candidate-parent searches and suitable spatial data structures, the observed cost can be substantially lower, closer to near-linear or mildly superlinear behavior in practical computations.

Thus the MSR algorithm is computationally feasible because exact topology optimization is confined to bounded-size subproblems, while the large-scale problem is treated by hierarchical initialization and energy-decreasing local and global refinements. Nevertheless, the method remains heuristic for large point sets. The algorithm does not provide a certificate of global optimality, and its output should therefore be interpreted as a high-quality numerical approximation to an $\alpha$-optimal transport path, rather than a certified global minimizer.

The transport paths in the following example were generated using the MSR algorithm.
\begin{example}
Let $\left\{ y_{i}\right\} $ be 50 random points in the square $%
\left[ 0,1\right] \times \left[ 0,1\right] $. Then, $\left\{ y_{i}\right\} $
determines an atomic probability measure
\begin{equation*}
\mathbf{b}=\sum_{i=1}^{50}\frac{1}{50}\delta _{y_{i}}.
\end{equation*}%
Let $\mathbf{a}=\delta _{O}$ where $O=\left( 0,0\right) $ is the origin.
Then an optimal transport path from $\mathbf{a}$ to $\mathbf{b}$ looks like
the following figures with $\alpha =1,0.75, 0.5$ and $0.25$ respectively:

\noindent
\includegraphics[height=1.25in]{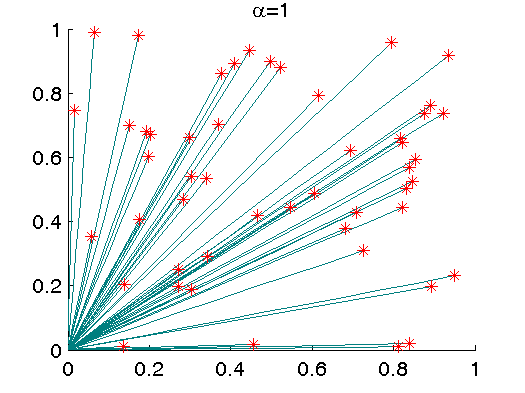}%
\includegraphics[height=1.25in]{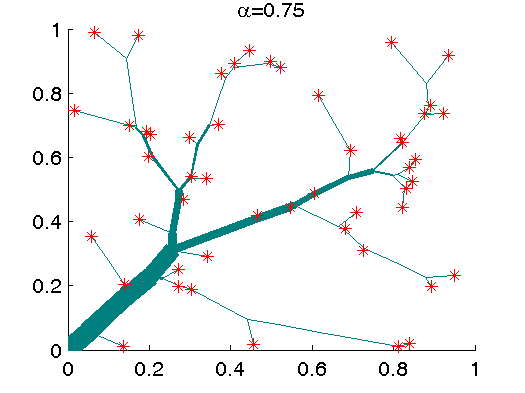}
\includegraphics[height=1.25in]{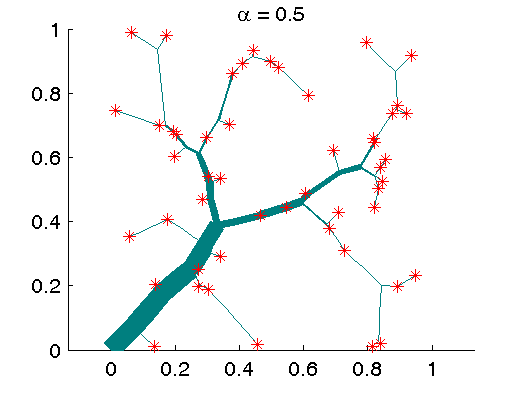}%
\includegraphics[height=1.25in]{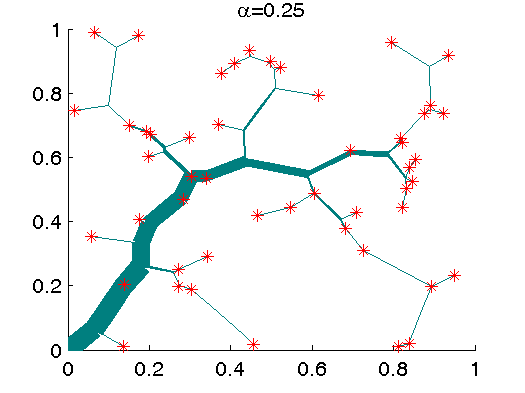}
\end{example}

\begin{example} Let $\mathbf{a}=m_{1}\delta _{x_{1}}+m_{2}\delta _{x_{2}}$ and $\mathbf{b}=m_{3}\delta_{x_{3}} $ with $m_{3}=m_{1}+m_{2}$. Then, in the non-degenerate case, the optimal transport path from $\mathbf{a}$ to $\mathbf{b}$ under the $\mathbf{M}_{\alpha }$ cost looks like the \textquotedblleft Y shaped\textquotedblright\ graph. 
\begin{center} \includegraphics[width=0.5\textwidth,height=2.5in]{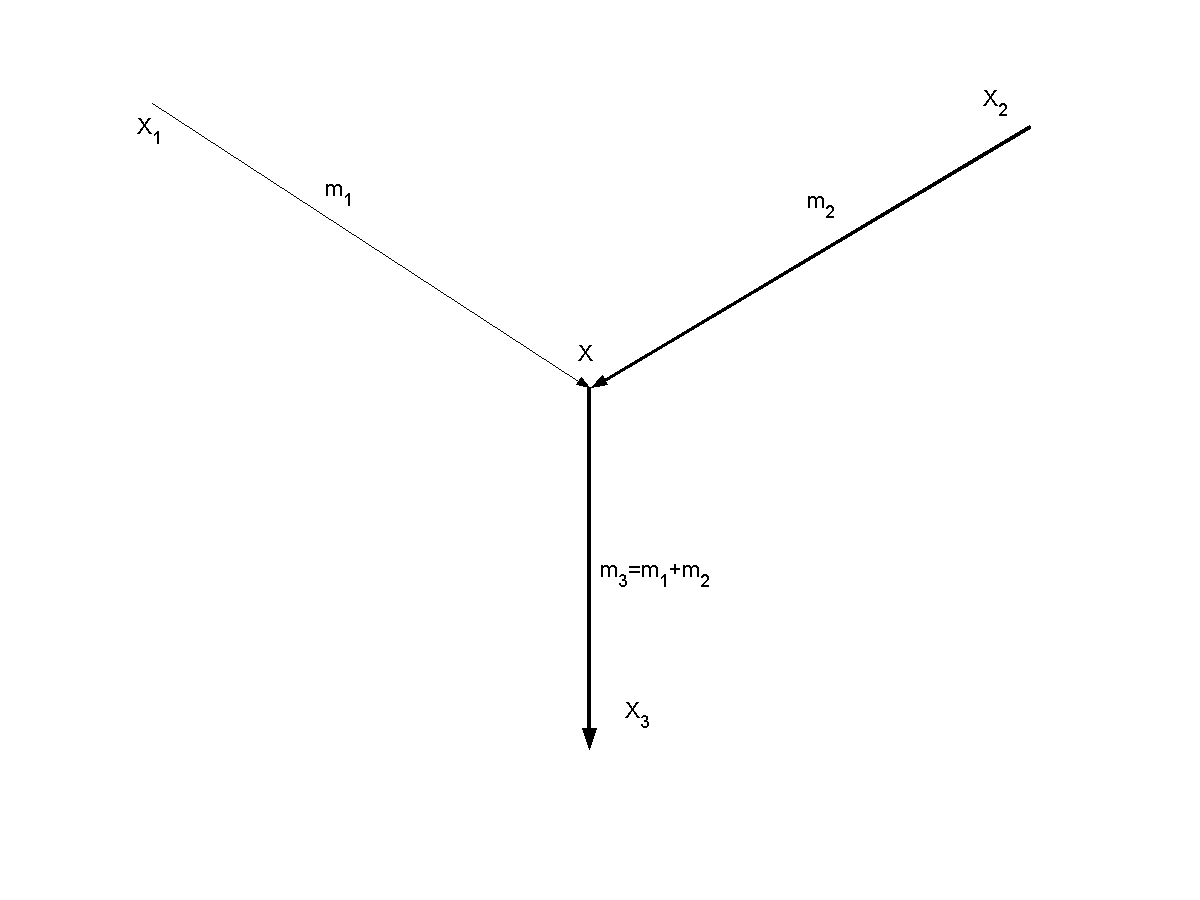}
\end{center}
\vspace{-10pt}
Here the interior vertex $x$ is determined by a balance formula: 
\begin{equation}
m_{1}^{\alpha }\overrightarrow{n_{1}}+m_{2}^{\alpha }\overrightarrow{n_{2}}%
+m_{3}^{\alpha }\overrightarrow{n_{3}}=\vec{0},  \label{balance_three}
\end{equation}%
where $\overrightarrow{n_{i}}=\frac{x_{i}-x}{\left\vert x-x_{i}\right\vert }$
is the unit vector from $x$ to $x_{i}$, $i=1,2,3$. Let $\theta _{i}$ be the
angle between $\overrightarrow{n_{i}}$ and $-$ $\overrightarrow{n_{3}}$ for $%
i=1,2$ and $k_{1}=\frac{m_{1}}{m_{1}+m_{2}}$, $k_{2}=\frac{m_{2}}{m_{1}+m_{2}%
}=1-k_{1}$.  Then, it is easy to find that the angles satisfy

\begin{equation*}
\cos \theta _{1}=\frac{k_{1}^{2\alpha }+1-k_{2}^{2\alpha }}{2k_{1}^{\alpha }}%
\text{ and }\cos \theta _{2}=\frac{k_{2}^{2\alpha }+1-k_{1}^{2\alpha }}{%
2k_{2}^{\alpha }}.
\end{equation*}%

Moreover,
\begin{equation}
\label{eqn: cos_sum}
\cos \left( \theta _{1}+\theta _{2}\right) =\frac{1-k_{1}^{2\alpha
}-k_{2}^{2\alpha }}{2k_{1}^{\alpha }k_{2}^{\alpha }}.
\end{equation}%
In particular, if $m_{1}=m_{2}$, then $k_1=k_2=\frac{1}{2}$ and
\begin{equation*}
\theta _{1}+\theta _{2}=\arccos \left( 2^{2\alpha -1}-1\right) .
\end{equation*}%
Note that for any $m_{1}>0$ and $m_{2}>0$, by (\ref{eqn: cos_sum}),
\begin{itemize}
\item if $\alpha =0$, then $\theta _{1}+\theta _{2}=120^o$;
\item if $\alpha =1/2$, then $\theta _{1}+\theta _{2}=90^o$;
\item if $\alpha =2/3$, then $\theta _{1}+\theta _{2}$ is nearly $75^o$.
\item Assume $\alpha>0$. If the ratio $\frac{m_2}{m_1}\rightarrow 0$, then $k_1\rightarrow 1, k_2\rightarrow 0$. In this case, $\theta_1\rightarrow 0$ indicates that the main stream nearly preserves its orientation.
\end{itemize}
\end{example}

In the more general case, we have the following results for angles between edges:
\begin{proposition}\cite[Proposition 2.1]{xia2011boundary}
 Let $G\in Path_{0}\left( \mathbf{a,b}\right) $ be any $\alpha $%
-optimal transport path with $\alpha \in \lbrack 0,1)$. Let $v\in V\left(
G\right) $ be an interior vertex of $G$, i.e. not one of the boundary
vertices $\left\{ x_{1},\cdots ,x_{k},y_{1},\cdots ,y_{\ell }\right\} $. Let 
$\left\{ e_{i}\right\} _{i=1}^{\deg \left( v\right) }$ be the edges in $G$
with $v$ being one of its endpoints. Let $w\left( e_{i}\right) $ be the
corresponding weight on $e_{i}$, and $\vec{e}_{i}$ be the unit directional
vector of the edge $e_{i}$ from $v$ to the other endpoint of $e_{i}$. Then,

\begin{enumerate}
\item There exists a balance equation at $v:$%
\begin{equation}
\sum_{i=1}^{\deg \left( v\right) }\left[ w\left( e_{i}\right) \right]
^{\alpha }\vec{e}_{i}=\vec{0}\text{.}
\label{balance_v}
\end{equation}

\item The minimum angle between any two edges in $\left\{ e_{i}\right\}
_{i=1}^{\deg \left( v\right) }$ is uniformly bounded below by 
\begin{equation}
\label{theta_alpha}
\theta _{\alpha }:=\left\{ 
\begin{array}{ll}
\frac{\pi }{2}, & \text{if }0<\alpha \leq \frac{1}{2}. \\ 
\arccos \left( 2^{2\alpha -1}-1\right) , & \text{if }\frac{1}{2}<\alpha <1%
\text{ or }\alpha =0%
\end{array}%
\right. .
\end{equation}
\item The degree $\deg \left( v\right) $ of $v$ is bounded above by some
constant 
\begin{equation}
D\left( \alpha ,m\right),
\label{D_alpha_m}
\end{equation}
 depending only on $\alpha $ and the
dimension $m$ of $\mathbb{R}^{m}$.
\end{enumerate}
\end{proposition}

For each $G \in Path_0 \left(\mathbf{ a,b}\right) $, we have the following trivial
but important lemma.

\begin{lemma}\cite[Lemma 2.1]{xia2003}
\label{lemma: mass_comparision}
Suppose $G\in Path_{0}\left(\mathbf{ a,b}\right) $ with $\mathbf{a,b}\in \mathcal{A}_{\Lambda }$ as in (\ref{a_lambda}).
Then for any edge $e\in E\left( G\right) $, we have 
$
0<w\left( e\right) \leq \Lambda .
$
Moreover,
$
\frac{\mathbf{M}_{\alpha }\left( G\right) }{\Lambda ^{\alpha }}\geq \frac{%
\mathbf{M}_1\left( G\right) }{\Lambda }
$ when $0\le \alpha\le 1$.
\end{lemma}

The following proposition provides an interesting formula of $\mathbf{M}_{\alpha
}\left( G\right) $ of an optimal transport path $G$ in terms of boundary
values. 
\begin{proposition}\cite[Proposition 6.4]{xia2014landscape}
For any $G\in Path\left( \mathbf{a,b}\right) $, it holds that 
\begin{equation}
\mathbf{M}_{\alpha }\left( G\right) =\sum_{v\in V\left( G\right) }\vec{%
\mathbf{m}}_{\alpha }\left( v\right) \cdot \mathbf{v,}  \label{cost_momentum}
\end{equation}
where $\mathbf{v}$ is the corresponding position vector of $v$ in $\mathbb{R}%
^{m}$ and 
\begin{equation*}
\vec{\mathbf{m}}_{\alpha }\left( v\right) =\sum_{e^{+}=v}w\left( e\right)
^{\alpha }\vec{e}-\sum_{e^{-}=v}w\left( e\right) ^{\alpha }\vec{e}. 
\end{equation*}%

If $G$ is an $\alpha $-optimal transport path, then a reformulation of the balance equation (\ref{balance_v}) gives $\vec{\mathbf{m}}_{\alpha }\left( v\right) =0$ for any interior vertex $v\in
V\left( G\right) \backslash \left\{ x_{1},\cdots ,x_{k},y_{1},\cdots
,y_{\ell }\right\} $ and thus 
\begin{equation}
\mathbf{M}_{\alpha }\left( G\right) =\sum_{v\in \left\{ x_{1},x_{2},\cdots
,x_{k},y_{1},\cdots ,y_{\ell }\right\} }\vec{\mathbf{m}}_{\alpha }\left(
v\right) \cdot \mathbf{v}.  \label{cost_momentum_boundary}
\end{equation}
\end{proposition}
Using (\ref{cost_momentum_boundary}) as well as the property that $\mathbf{M}_{\alpha}(G)$ is translational invariant, one gets
a necessary condition for a transport path being optimal:
\begin{corollary}\cite[Corollary 6.5]{xia2014landscape}
Suppose $G$ is an $\alpha$-optimal transport path from $\mathbf{a}$ to $%
\mathbf{b}$. Then it holds that 
\begin{equation*}
\sum_{v\in \left\{ x_{1},x_{2},\cdots ,x_{k},y_{1},\cdots ,y_{\ell }\right\}
}\vec{\mathbf{m}}_{\alpha }\left( v\right) =0. 
\end{equation*}
\end{corollary}

\subsection{
Continuous Case: General Transport paths}

We aim to generalize the notion of a transport path from atomic measures to arbitrary measures, capturing branching structures in a more general setting. We first express general transport paths in terms of vector measures.
\subsubsection{Set-up and Existence result in terms of vector measures \cite{xia2003}}

Let $\mathcal{M}_\Lambda(X)$ be the space of Radon measures $\mu$ on $X$ with total mass $\mu(X)=\Lambda$. When $\Lambda=1$, we also denote $\mathcal{M}_1(X)$ by $\mathcal{P}(X)$, the space of all probability measures on $X$.

Recall that in Remark \ref{remark:path}, each transport path between atomic measures can be viewed as a vector-valued measure on $X$. This enables us to consider transport paths between two Radon measures as follows:  

Let $\mu$, $\nu\in \mathcal{M}%
_{\Lambda }\left( X\right) $ be any two Radon measures on $X$ with equal
total mass $\Lambda $.

\begin{definition}
A \emph{transport path} from $\mu$ to $\nu$ is a vector-valued measure $T$ on $X$ 
such that
\[
\mathrm{div}(T) = \mu - \nu
\]
in the sense of distributions.
\end{definition}
 \noindent This expresses the conservation of mass at every point in $X$, generalizing the Kirchhoff condition from discrete vertices to a continuous domain.
 
 Now, we use a limit process to define the $\mathbf{M}_\alpha$ cost of a general transport path as follows. Given a vector-valued measure $T\in \mathcal{M}^{m}\left( X\right) $, \textbf{an approximating graph sequence for }$T$ is a sequence of triples $\left\{ \mathbf{a}_{i}, \mathbf{b}_{i}, G_{i}\right\} $ such that 
\begin{itemize}
    \item $\left\{ \mathbf{a}_{i}\right\} $ and $\left\{ \mathbf{b}_{i}\right\} $ are
two sequences of atomic measures in $\mathcal{A}_{\Lambda }\left( X\right) $ with $\mathbf{a}_{i}\rightharpoonup \mu,\ \ \mathbf{b}_{i}\rightharpoonup \nu$  weakly as Radon measures;
\item $G_{i} \in Path\left(
\mathbf{a}_{i}, \mathbf{b}_{i}\right) $ with  $G_{i}\rightharpoonup T$ weakly as vector-valued measures.
\end{itemize}

Note that such approximating graph sequences are easily found by, for instance, using dyadic approximation of Radon measures with atomic measures, and corresponding transport paths between atomic measures.

Now, given any $\alpha \in \left[ 0,1\right] $,
we define its $\mathbf{M}_{\alpha }$ cost to be
\begin{equation*}
\mathbf{M}_{\alpha }\left( T\right) :=\inf \{\lim \inf_{i\rightarrow \infty }%
\mathbf{M}_{\alpha }\left( G_{i}\right)\} ,
\end{equation*}%
where the infimum is over the set of all possible approximating graph
sequence $\left\{ a_{i},b_{i},G_{i}\right\} $ of $T$.

\begin{example}
\label{ex:dyadic-path}
Suppose $\mu \in \mathcal{M}_{\Lambda}\left( X\right) $ is supported in a cube $Q\subset \mathbb{R}^m$ whose center is denoted by $c_Q$.  By using a sequence of dyadic approximations of $\mu$ as shown in Figure \ref{fig:dyadic}, one can construct a transport path $T$ from $\Lambda\delta_{c_Q}$ to $\mu$ such that for any $\alpha \in (1-\frac{1}{m},1]$,
\begin{equation*}
\mathbf{M}_{\alpha }\left( T\right) \leq \frac{\Lambda^\alpha}{2^{1-m\left( 1-\alpha
\right) }-1}\frac{diam(Q)}{2}.
\end{equation*}%
\begin{figure}
    \centering
    \includegraphics[width=0.5\linewidth]{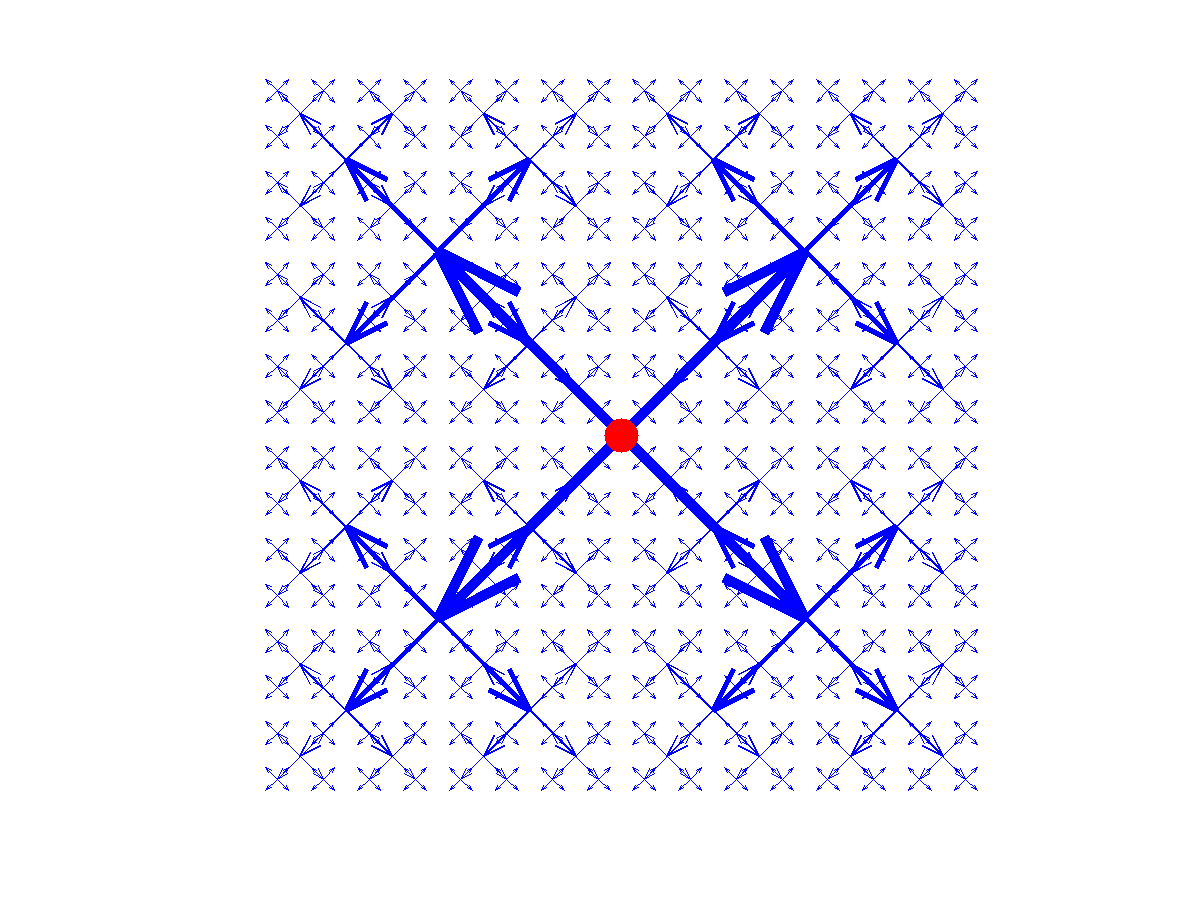}
    \caption{Using dyadic approximation to construct a transport path from a Dirac mass at the center of a cube to the Lebesgue measure on the cube.}
    \label{fig:dyadic}
\end{figure}
\end{example}

\begin{remark}
\label{remark: dimension}
The construction and estimation given in the example works for all Radon measures. For a particular measure, one may find other constructions with finite  $\mathbf{M}_\alpha$ cost for lower values of $\alpha$. The critical value of all possible such $\alpha$ depends on the transport dimension of the measure. We refer to \cite{xia2010dimension, xia2011metric} for details of the transport dimension.
\end{remark}

\subsubsection{In terms of rectifiable flat 1-chains}

To get a better description of the geometric structures of those general transport paths, we express them in the language of geometric measure theory, in particular, in rectifiable currents and flat chains. 

We first recall some terminology about rectifiable currents as in \cite{Federer1969} or \cite{simon}.

Let $\Omega \subseteq \mathbb{R}^{m}$ be an open subset and $\mathcal{D}%
^{k}(\Omega )$ be the set of all $C^{\infty }$ differential $k$-forms in $%
\Omega $ with compact support with the usual Fr\'{e}chet topology \cite{Federer1969}%
. An $k$\textit{-dimensional current }$S$\textit{\ in }$\Omega $\textit{\ }%
is a continuous linear functional on $\mathcal{D}^{k}(\Omega )$. Let $%
\mathcal{D}_{k}(\Omega )$ denote the set of all $k$-dimensional currents in $%
\Omega $. Motivated by Stokes' theorem, the \textit{boundary} of a current $%
S\in \mathcal{D}_{k}(\Omega )$ is the $k-1$ dimensional current defined by
$
\partial S\left( \psi \right) :=S\left( d\psi \right)
$
for any $C^{\infty }$ form $\psi \in \mathcal{D}%
^{k-1}(\Omega )$. A sequence of currents $S_{i}\in \mathcal{D}_{k}(\Omega )$ is said
to be weakly convergent to another current $S\in \mathcal{D}_{k}(\Omega )$,
denoted by $S_{i}\rightharpoonup S$, if
$
S_{i}(\psi )\rightarrow S(\psi )
$
for any $\psi \in \mathcal{D}^{k}(\Omega )$.

As in \cite{simon}, a subset $M\subseteq \mathbb{R}^{m}$ is called
(countably) $k-$\textit{rectifiable} if $M=\bigcup\limits_{i=0}^{\infty
}M_{i}$, where $\mathcal{H}^{k}\left( M_{0}\right) =0$ under the $k-$%
dimensional Hausdorff measure $\mathcal{H}^{k}$ and each $M_{i}$, for $%
i=1,2,\cdots ,$ is a subset of an $k-$dimensional $C^{1}$ submanifold in $%
\mathbb{R}^{m}$. A \textit{rectifiable current} $S$ is a current coming from
an oriented rectifiable set with multiplicities. More precisely, we have
\begin{definition}[Rectifiable Currents]
A $k$-dimensional current $T$ in $\Omega$ is \emph{rectifiable} if it can be written in the form
\[
T(\omega) = \int_{M} \langle \omega(x), \xi(x) \rangle \,\theta(x)\, d\mathcal{H}^k(x), \forall
\omega \in \mathcal{D}^{k}(\Omega )
\]
where:
\begin{itemize}
    \item $M\subset \Omega$ is countably $k$-rectifiable,
    \item $\theta $ is a $\mathcal{H}^{1}\lfloor M$ integrable positive
function, called the multiplicity function of $T$,
    \item $\xi(x)$ is a unit simple $k$-vector giving orientation of $T$.
\end{itemize}
\end{definition}
The rectifiable current $T$ described as above is often denoted by
$
T=\underset{=}{\tau }(M,\theta ,\xi )$. 
The mass of $T$ is $\mathbf{M}(T) = \int_{M} \theta(x)\, d\mathcal{H}^k(x)$. For instance, each weighted directed graph $G\in Path\left( \mathbf{a,b}%
\right) $ in the discrete case determines a rectifiable 1-current $G=\underset{=}{\tau }\left(
\mathbf{G},w,\xi \right) $ in $\mathbb{R}^{m},$ with $\partial G=\mathbf{b-a}
$ as currents.

\begin{definition}[Flat Chains]
Let $k\in\mathbb{N}$. A \emph{$k$-dimensional flat chain} $T$ in $\mathbb{R}^d$ is a current that can be approximated in the flat norm
\[
\mathcal{F}(T) := \inf\{ \mathbf{M}(R) + \mathbf{M}(S) \,:\, T = R + \partial S, \, R \text{ a $k$-current},\, S \text{ a $(k+1)$-current} \}
\]
by polyhedral $k$-chains with real coefficients.  
Here $\mathbf{M}(\cdot)$ denotes the mass norm.  
\end{definition}

It is well known that a real flat chain with finite mass may not be rectifiable. Nevertheless, the following result holds: 

\begin{theorem}[Rectifiable Slicing Theorem \cite{White1999}]
\label{rectifiableslices}\bigskip Let $T$ be a $k$-flat chain of finite mass
in $\mathbb{R}^m$. The following are equivalent:
\begin{enumerate}
\item $T$ is rectifiable.

\item For almost every $\left( m-k\right) $-plane $P$ parallel to a
coordinate plane, the slice $T\cap P$ is a rectifiable $0$-chain.
\end{enumerate}
\end{theorem}

We now aim at expressing transport paths as rectifiable flat 1-chains. Through a reference from Thierry De Pauw, the author learned of a result of Brian White~\cite{White1999}: every real flat $0$-chain with finite $\mathbf M_\alpha$ mass, for $\alpha\in[0,1)$, is rectifiable. This leads to the result:
\begin{proposition}[\cite{xia2004interior}]
\label{prop:rectfiability}
    Let $T$ be a transport path with finite $\mathbf{M}_\alpha$ cost for some $\alpha\in[0,1)$. Then, $T$ is a rectifiable flat 1-chain.
\end{proposition}

\begin{proof}[Sketch of the proof:] Using Brian White's result and the rectifiable slicing theorem, one can show that any real flat $k$-chain $T$ of finite $\mathbf{M}$
mass and finite $\mathbf{M}_{\alpha }$ mass is rectifiable. Now, let $T$ be a transport path from $\mu$ to $\nu$ with finite $\mathbf{M}_\alpha$ mass as described above in terms of vector measures. Using the equivalence between weak convergence and flat norm convergence \cite{simon}, it follows that $T$ is a real flat 1-chain with $\partial T=\nu -\mu$. 
By means of Lemma \ref{lemma: mass_comparision}, it follows that finite $\mathbf{M}_\alpha$ 
mass also implies finite $\mathbf{M}$ mass. As a result, $T$ is rectifiable.
\end{proof}

By the proposition, $T$ can be expressed as a rectifiable flat 1-chain
\begin{equation*}
T=\underset{=}{\tau }(M,\theta ,\xi )
\end{equation*}
whose $\alpha$-mass is expressed as
\[
\mathbf{M}_\alpha(T) = \int_{M} \theta(x)^\alpha\, d\mathcal{H}^1(x).
\]
 This generalizes the discrete sum in (\ref{alpha_cost}) to the continuum setting.

Let
\begin{equation*}
Path\left( \mu,\nu\right)=\{T \text{ is a real rectifiable flat 1-chain with } \partial T=\nu -\mu \}
\end{equation*}
be the space of all transport paths from $\mu$ to $\nu$ with finite $\mathbf{M}_\alpha$-cost.

\begin{problem}[Ramified optimal transport problem: Continuous Version]  Given two Radon measures $\mu,$
$\nu\in \mathcal{M}_{\Lambda }\left( X\right) $ on $X\subset \mathbb{R}%
^{m}$ and $\alpha \in [0,1)$, minimize $\mathbf{M}_{\alpha}(T)$ among all $T\in Path(\mu, \nu)$.
\label{problem:continuous}
\end{problem}
In terms of Plateau-type problem, it is: Minimize $\mathbf{M}_{\alpha}(T)$ among all real rectifiable flat 1-chains with $\partial T=\nu-\mu$.

Now, we restate \cite[Theorem 3.1]{xia2003} as follows.
\begin{theorem}[Existence of $\alpha$-Optimal Transport Paths]
Let $\mu$ and $\nu$ be two positive Radon measures on a convex compact set $X \subset \mathbb{R}^m$ with equal total mass $\Lambda > 0$. For any $\alpha \in (1-\frac{1}{m},1]$, there exists a transport path $S$ from $\mu$ to $\nu$ that minimizes the $\mathbf{M}_\alpha$ cost among all transport paths, i.e.,
\[
\mathbf{M}_\alpha(S) = \inf\{ \mathbf{M}_\alpha(T) \, : \, \text{$T$ is a transport path from $\mu$ to $\nu$}\}.
\]
Moreover,
\begin{equation*}
\mathbf{M}_{\alpha }\left( S\right) \leq \frac{\Lambda ^{\alpha }}{%
2^{1-m\left( 1-\alpha \right) }-1}\frac{\sqrt{m}d}{2},
\end{equation*}
where $d$ is the diameter of the convex hull of supports of $\mu$ and $\nu$.
\end{theorem}

\begin{proof}
We sketch a standard proof using the ``direct method in the calculus of variations''.  

\textit{Step 1: Minimizing sequence.} Let $\{T_n\}_{n \in \mathbb{N}}$ be a sequence of transport paths such that
\[
\mathbf{M}_\alpha(T_n) \to \inf_{T\in Path(\mu,\nu)} \mathbf{M}_\alpha(T) < \infty.
\]
where the boundedness comes from Example \ref{ex:dyadic-path}.

\textit{Step 2: Compactness.} Since $\{T_n\}$ is uniformly bounded in $\mathbf{M}_\alpha$-mass, by Lemma \ref{lemma: mass_comparision}, it is also  uniformly bounded in $\mathbf{M}$-mass. Since $\{T_n\}$
is supported in the compact set $X$, by standard compactness results for 1-dimensional rectifiable currents or vector-valued measures or flat 1-chains (see, e.g., Federer \cite{Federer1969}), there exists a subsequence (still denoted $T_n$) that converges weakly-* to a limit $S$.

\textit{Step 3: Boundary condition.} The divergence operator (or the boundary operator) is continuous under weak-* convergence of measures. Hence, 
\[
\mathrm{div}(S) = \lim_{n \to \infty} \mathrm{div}(T_n) = \lim_{n \to \infty} (\mu - \nu) = \mu - \nu.
\]
Thus $S$ is a transport path from $\mu$ to $\nu$.

\textit{Step 4: Lower semicontinuity.}
\textcolor{blue}{
By the lower semicontinuity of the relaxed $\mathbf M_\alpha$ functional
with respect to weak convergence of currents
\cite{ColomboDeRosaMarcheseStuvard2017}, we obtain
\[
\mathbf M_\alpha(S)
\leq
\liminf_{n\to\infty}\mathbf M_\alpha(T_n).
\]
Since $\{T_n\}$ is a minimizing sequence, it follows that
\[
\mathbf M_\alpha(S)
\leq
\liminf_{n\to\infty}\mathbf M_\alpha(T_n)
=
\inf_T \mathbf M_\alpha(T).
\]
Thus $S$ is an $\alpha$-optimal transport path from $\mu$ to $\nu$.
 The upper bound for $\mathbf M_\alpha(S)$ follows from Example~1.17.}
\end{proof}

\begin{remark}
As seen from the above proof, for any $\alpha\in [0,1]$, there exists an $\alpha$-optimal transport path from $\mu$ to $\nu$, provided that there exists some transport path with finite $\mathbf{M}_{\alpha}$-cost. The latter information, as stated in Remark \ref{remark: dimension} in general, depends on the transport dimension of the signed measure $\nu-\mu$. In the above proof, we limit $\alpha$ to within $(1-\frac{1}{m}, 1]$ because Example \ref{ex:dyadic-path} always provides a finite $\alpha$-cost transport path.
\end{remark}

\subsubsection{Regularity of optimal transport paths}
After studying the existence of $\alpha$-optimal transport paths, one would like to know their regularity, i.e. to see how nice these paths are. Under the assumption that one of $\mu$ and $\nu$ is atomic, we have the following interior regularity theorem \cite{xia2004interior}:
\begin{theorem}
\label{thm:interior}
Suppose $T\in Path\left( \mu,\nu\right) $ is an optimal transport
path with finite $\mathbf{M}_{\alpha }$ cost. For any point $\xi \in
spt\left( T\right) \setminus spt\left( \mu\cup \nu\right) $, there
is an open neighborhood $B_{\xi }$ of $\xi $, such that
$
T\lfloor B_{\xi }
$
is a cone that consists of finitely many line segments with suitable multiplicities.
\end{theorem}

The following corollary implies that our approach in the continuous setting is a generalization of the one in the discrete setting. 
\begin{corollary}[Discrete Minimizers are Finite Weighted Graphs of Line Segments]
\label{cor:discrete-graph}
Let $\alpha \in [0,1)$ and let $\mathbf{a},\mathbf{b} \in \mathcal{A}_\Lambda(X)$ be atomic measures.  
Any $\alpha$-optimal transport path $T$ from $\mathbf{a}$ to $\mathbf{b}$ is supported on a finite $1$-rectifiable set that is a finite union of straight line segments. More precisely:
\begin{enumerate}
    \item $\mathrm{spt}(T)$ is a finite weighted directed graph $G=(V,E,w)$ with
    \[
    \{x_1,\dots,x_k\} \cup \{y_1,\dots,y_\ell\} \subset V,
    \]
    each edge $e\in E$ being a straight segment of $\mathbb{R}^d$.
    \item $G$ is acyclic (hence a forest).
    \item (Kirchhoff balance) At every vertex $v\in V$,
    \[
    \sum_{\substack{e \in E\\ e^- = v}} w(e)
    \;=\;
    \sum_{\substack{e \in E\\ e^+ = v}} w(e)
    \;+\;
    \sum_{\substack{v = x_i}} m_i
    \;-\;
    \sum_{\substack{v = y_j}} n_j.
    \]
    \item (\emph{First-order optimality at junctions})  
    At each interior vertex $v$,
    \[
    \sum_{e:\ e \text{ is an edge with one endpoint } v} w(e)^{\alpha}\,\vec{n}_e(v) \;=\; 0,
    \]
    where $\vec{n}_e(v)$ is the unit directional vector of  $e$ at $v$, oriented outward from $v$.
\end{enumerate}
\end{corollary}
\begin{proof}[Proof Sketch]
By Proposition \ref{prop:rectfiability} and Theorem \ref{thm:interior}, any minimizer $T$ is $1$-rectifiable with finite $\mathcal{H}^1$-measure and finitely many junctions.  
Stationarity implies edges are straight; cycles are eliminated by cost reduction for $\alpha<1$.  
The balance law follows from $\partial T=\mathbf{b}-\mathbf{a}$.  
The vertex equilibrium condition follows from first-order variations.
\end{proof}
In general, the support of an optimal transport path $T$
may not necessarily be $1$ dimensional nearby its boundary, for the support
of the measure $\mu-\nu$ may even contain an open subset of $%
\mathbb{R}^m$. For instance, one may take $\mu$ to be some Lebesgue
measure on a domain $\Omega $ and $\nu$ to be some atomic measure on $%
\Omega $. Then, the support of $\mu-\nu$ has the same dimension of
$\Omega $, which is not necessarily $1$ dimensional. So, the question is how
to describe the behavior of $T$ when a carrying set of $T$ is possibly dense
in the whole space $X$.

To study the boundary regularity of an optimal transport path, the main idea
of our approach is to study its superlevel sets. This idea is motivated from
observing vein structure of a tree leaf provided by the nature. Here,
we show that each superlevel set of an optimal transport path is
locally concentrated on the support of an $1$-dimensional bi-Lipschitz
chain, which is analogous to vein structures of a tree leaf.

We first clarify some terminology. For any $\lambda >0$, the $\lambda$-\textit{superlevel set} of a rectifiable current $T=\underset{=}{%
\tau }\left( M,\theta ,\xi \right) $ is the set
\begin{equation}
T_{\lambda }=\left\{ x\in M:\theta \left( x\right) \geq \lambda \right\} .
\label{levelset}
\end{equation}

Also, a \textit{bi-Lipschitz chain} $P$ is a finite sum of bi-Lipschitz
curves in $X$ with real coefficient multiplicities. That is,
\begin{equation}
P=\sum\limits_{i=1}^{K}m_{i}\Gamma _{i}  \label{bilipschitz}
\end{equation}%
for some real numbers $m_{i}>0$, and some bi-Lipschitz curves $\Gamma _{i}$
with $i=1,2,\cdots ,K$ in $X$. The support $spt\left( P\right) $ of the
bi-Lipschitz chain $P$ is the union of the image of every bi-Lipschitz curve
$\Gamma _{i}$ in $X$.

Now, we state our boundary regularity theorem as follows:

\begin{theorem}\cite[Theorem 4.1]{xia2011boundary}
\label{theorem1}For any $\mu,\nu\in \mathcal{M}_{\Lambda }\left(
X\right) $, let $T=\underset{=}{\tau }\left( M,\theta ,\xi \right) \in
Path\left( \mu,\nu\right) $ be any $\mathbf{\alpha -}$optimal
transport path with finite $\mathbf{M}_{\alpha }$ cost for some $0\leq
\alpha <1$. Then, for any $\epsilon >0$ and any $p\in M$, there exists an
open ball neighborhood $B_{r}\left( p\right) $ about $p$ and a decomposition
\begin{equation*}
T \hspace{.1\baselineskip} \rule{.5pt}{.5 \baselineskip}\rule%
{.5\baselineskip}{.5pt} \hspace{.2\baselineskip} B_{r}\left( p\right) =P+R
\end{equation*}%
as $1-$dimensional rectifiable currents such that

\begin{itemize}
\item[a)] $P$ is a bi-Lipschitz chain in the form of (\ref{bilipschitz}).

\item[b)] $R\in Path\left( \mu _{R},\nu _{R}\right) $ is an $1-$%
dimensional rectifiable current for some Radon measures $\mu _{R}$ and $%
\mu _{R}$ with mass $\left| \left| \mu _{R}\right| \right| =\left|
\left| \mu _{R}\right| \right| <\epsilon $.

\item[c)] Moreover, inside the ball $B_{r}\left( p\right) $, the $\epsilon -$
superlevel set $T_{\epsilon }$ of $T$ as defined in (\ref{levelset}) is a
subset of the support of the bi-Lipschitz chain $P$.
\end{itemize}
\end{theorem}
Here, $P=P_\epsilon$ plays a similar role as in the vein structure of a tree leaf. As $\epsilon>0$
decreases, one can see more details about the structure.

\subsubsection{The $d_\alpha$ metric}
At the end of this section, we want to describe the geometry of the space of probability measures under the metric derived from ramified optimal transportation.

\begin{definition}
    For any $\alpha \in [0,1]$, we define
\begin{equation*}
d_{\alpha }\left( \mu,\nu\right) :=\inf \left\{ \mathbf{M}_{\alpha
}\left( T\right) :T\in Path\left( \mu,\nu\right) \right\} ,
\end{equation*}%
for any two Radon measures $\mu,$ $\nu\in \mathcal{M}_{\Lambda
}\left( X\right)$.
\end{definition}
Note that, if there exists an $\alpha$-optimal transport path $S$ from $\mu$ to $\nu$. Then,
\begin{equation*}
d_{\alpha }\left( \mu,\nu\right) =\mathbf{M}_\alpha(S)=\min \left\{ \mathbf{M}_{\alpha
}\left( T\right) :T\in Path\left( \mu,\nu\right) \right\}.
\end{equation*}%
This is the case in particular if $\alpha>1-\frac{1}{m}$.

Also note that for any $\Lambda >0$ and any $\mu$, $\nu\in \mathcal{M}%
_{\Lambda }\left( X\right) $,
\begin{equation*}
d_{\alpha }\left( \mu,\nu\right) =\Lambda ^{\alpha }d_{\alpha
}\left( \frac{\mu}{\Lambda },\frac{\nu}{\Lambda }\right) .
\end{equation*}
Thus, we may assume $\Lambda =1$. Some key results about $d_{\alpha}$ is stated in the following theorem:

\begin{theorem}\cite[Theorems 4.1, 4.2, 5.1]{xia2003}
$d_{\alpha }$ is a metric on $\mathcal{M}_{1}\left( X\right)$ and metrizes the weak * topology of $\mathcal{M}_{1}\left(
X\right)$. Moreover, the space $\left( \mathcal{M}_{1}\left( X\right), d_{\alpha }\right) $ is a length space in the sense that for any $%
\mu,\nu\in \mathcal{M}_{1}\left( X\right) $, each $\alpha$-optimal transport path $T$ corresponds to a
continuous map
\begin{equation*}
\psi :\left[ 0,d_{\alpha }\left( \mu,\nu\right) \right]
\rightarrow \mathcal{M}_{1}\left( X\right)
\end{equation*}
such that $\psi \left( 0\right) =\mu,\psi \left( d_{\alpha }\left( \mu
^{+},\nu\right) \right) =\nu$ and for any $0\leq s_{1}<s_{2}\leq
d_{\alpha }\left( \mu,\nu\right) $,
\begin{equation*}
d_{\alpha }\left( \psi \left( s_{1}\right) ,\psi \left( s_{2}\right) \right)
=s_{2}-s_{1}\text{.}
\end{equation*}
\end{theorem}

How special is this metric? 
It turns out that this metric is an intrinsic metric induced by a quasi-metric. Let us describe it in details now.

A quasimetric space is a generalized metric space, in which
the distance satisfy a relaxed triangle inequality $d(x,y)\leq
c(d(x,z)+d(z,y))$ for some $c\geq 1$, rather than the usual triangle
inequality. Quasimetric concepts are critical for studying harmonic
analysis on spaces of homogeneous type. 
Further exploration of optimal transport paths provides us a
motivational example to study the geodesic problem in quasimetric
spaces in \cite{xia2009quasimetric}.  We first introduce a quasimetric on the space of atomic probability measures as follows: 
\begin{definition}
For any two atomic probability measures 

\begin{equation*}
\mathbf{a}=\sum_{i=1}^{m}a_{i}\delta _{x_{i}}\text{ and }\mathbf{b}%
=\sum_{j=1}^{n}b_{j}\delta _{y_{j}}
\end{equation*}
 on a metric space $(Y,d)$, and $%
\alpha <1$, define
\begin{equation*}
J_{\alpha }\left( \mathbf{a},\mathbf{b}\right) :=\min \left\{ 
\sum_{i=1}^{m}\sum_{j=1}^{n}\left( \gamma_{ij}\right) ^{\alpha }d\left( x_{i},y_{j}\right)
:\gamma =\sum_{i=1}^{m}\sum_{j=1}^{n}\gamma _{ij}\delta _{\left(
x_{i},y_{j}\right) }  \in Plan\left( \mathbf{a},\mathbf{b}\right)
\right\} .
\end{equation*}
\end{definition}

For any given natural number $N\in \mathbb{N}$ , let $\mathcal{A}_{N}(Y)$ be
the space of all atomic probability measures
$
\sum_{i=1}^{m}a_{i}\delta _{x_{i}}
$
on $Y$ with $m\leq N$.
One can show that $J_{\alpha }$ defines a quasimetric on $\mathcal{A}_{N}\left( Y\right) $. Note that, in general, $J_{\alpha }$ may fail to be a metric on $\mathcal{A}_{N}\left( Y\right) $ as demonstrated in the following example.

\begin{example}
For any $\alpha <1$, let $y$ be a positive real number. Then, we consider
three atomic measures in $Y=\mathbb{R}^{2}:$
\begin{equation*}
\mathbf{a}=\frac{1}{2}\delta _{\left( -1,y+1\right) }+\frac{1}{2}\delta
_{\left( 1,y+1\right) },\mathbf{b}=\delta _{\left( 0,0\right) }\text{ and }%
\mathbf{c}=\delta _{\left( 0,y\right) }\text{.}
\end{equation*}%
Then,
\begin{eqnarray*}
&&J_{\alpha }\left( \mathbf{a},\mathbf{c}\right) +J_{\alpha }\left( \mathbf{c%
},\mathbf{b}\right) -J_{\alpha }\left( \mathbf{a},\mathbf{b}\right) \\
&=&2\left( \frac{1}{2}\right) ^{\alpha }\sqrt{2}+y-2\left( \frac{1}{2}%
\right) ^{\alpha }\sqrt{1+\left( y+1\right) ^{2}}<0
\end{eqnarray*}%
whenever $y$ is large enough. Thus,\thinspace $J_{\alpha }$ does not satisfy the triangle inequality.
\end{example}

As we know, the intrinsic metric induced by a metric has
played an important role in the study of metric geometry. A natural
question is: will a quasimetric be able to induce an intrinsic
\textit{metric}?

To answer this question, we first studied properties of a quasimetric space $(X,J)$ in \cite{xia2009quasimetric}. Note that, as $J$ only satisfies the relaxed triangle inequality, $J$ is not necessarily continuous. One can easily show that $J_{\alpha}$ is not continuous by constructing a sequence of measures located on a "Y-shaped" optimal transport path.  Nevertheless, we showed that many well-known results in metric spaces still hold in suitable quasimetric
spaces. For instance,  the Ascoli-Arzel\`{a} theorem holds when $(X,J)$ is a complete quasimetric space and $J$ is lower semi-continuous.

Also, to study the geodesic problem in a quasimetric space $(X,J)$, one needs to define the length of a continuous curve $f: [a,b]\rightarrow (X,J)$. In a metric space $(X,d)$, we usually define 
\begin{equation*}
L\left( f\right) =\sup_{P}\sum_{i=1}^{N}d\left( f\left( t_{i-1}\right) ,f\left(
t_{i}\right) \right),
\end{equation*}%
where the supremum is over all partitions $P=\left\{
a=t_{0}<t_{1}<\cdots <t_{N}=b\right\} $ of $\left[ a,b\right] $.
However, on a quasimetric space, subdivision of a partition may {\it decrease} rather than increase the variation of $f$. One can not take the supremum of variations as before. Nevertheless, in the study of ramified optimal transportation, we observed in \cite[Remark 4.16]{xia2009quasimetric} that for each transport $G$ between two atomic probability measures in $\mathcal{A}_N(Y)$, there exists a piecewise metric Lipschitz curve $f$ in the quasimetric space $(\mathcal{A}_N(Y), J_{\alpha})$ such that the transport cost
   \[\mathbf{M}_{\alpha}(G)=\int_0^1 |\dot{f}(t)|_{J_{\alpha}}dt\]
   where the quasimetric derivative
\begin{equation*}
\left| \dot{f}\left( t\right) \right| _{J_{\alpha}}:=\lim_{s\rightarrow t}\frac{%
J_{\alpha}\left( f\left( s\right) ,f\left( t\right) \right) }{\left| s-t\right| }
\end{equation*}
provided the limit exists. This motivates us to define the length of a curve $f$ in a general quasimetric space $(X,J)$ to be
\begin{equation*}
L_{J}\left( f\right) :=\int_{a}^{b}\left| \dot{f}\left( t\right) \right|
_{J}dt
\end{equation*}
whenever it is well defined. 
  In particular,  for any natural number $M$,  every curve $f$ in the family $\mathcal{P}_M([a,b],(X,J))$ of all $M$-piecewise metric Lipschitz curves in $(X,J)$ has a well-defined length $L_J(f)$.  For any $x,y\in X$, we consider
the geodesic problem
\begin{equation}
\min \{ L_J\left( f\right) \}  \label{geodesic_problem}
\end{equation}%
among all $f$ in the family
\begin{equation*}
Path_{M}\left( x,y\right) =\left\{ f\in \mathcal{P}_{M}\left( \left[ 0,1%
\right] \text{,}\left( X,J\right) \right) \text{ with }f\left( 0\right)
=x;f\left( 1\right) =y\right\} .
\end{equation*}

\begin{definition}
Suppose $J$ is a quasimetric on $X$. For any $x,y\in X$, and $M\in \mathbb{N}
$, define
\begin{equation*}
D_{J}^{\left( M\right) }\left( x,y\right) =\inf \left\{ L_{J}\left( f\right)
:f\in Path_{M}\left( x,y\right) \right\}
\end{equation*}%
whenever $Path_{M}\left( x,y\right) $ is not empty, and set $D_{J}^{\left(
N\right) }\left( x,y\right) =\infty $ when $Path_{M}\left( x,y\right) $ is
empty. Since $D_{J}^{\left( M\right) }\left( x,y\right) $ is a decreasing
function of $M$, we define
\begin{equation*}
D_{J}\left( x,y\right) =\lim_{M\rightarrow \infty }D_{J}^{\left( M\right)
}\left( x,y\right) .
\end{equation*}
\end{definition}

For each quasimetric $J$ on $X$, $D_{J}$ is automatically a pseudometric.
When $D_{J}$ is indeed a metric on $X$, it is called the \textit{%
intrinsic metric}, or \textit{geodesic distance}, on $X$ induced by the
quasimetric $J$.

We explored some sufficient conditions on $X$ and $J$ that ensure $D_J$ to be a metric. One of them says that if $D_{J}^{\left( M\right)}(x,y)$ becomes a real valued constant $D_J(x,y)$ (i.e. independent of $M$) whenever $M$ is large enough, then $D_J$ is a metric on $X$. If in addition, the geodesic problem (\ref{geodesic_problem}) has a solution for $M$ large enough, then $(X,D_J)$ is a length space. It turns out the quasimetric space $(\mathcal{A}_N(Y), J_{\alpha})$ satisfies these conditions, and thus $J_{\alpha}$ induces a metric $D_{J_{\alpha}}$. The metric $d_{\alpha}$ coincides with this induced metric  $D_{J_{\alpha}}$ on $\mathcal{A}_N(Y)$. Consequently, an optimal transport path between two atomic
probability measures is simply a geodesic in the length space $(\mathcal{A}_N(Y), d_{\alpha})$  in the usual sense of metric
geometry.

\begin{remark}
Another question that we studied is the shape of the unit ball in the space of probability measures under the $d_\alpha$ metric. Using the landscape function to be described later, we showed in \cite{PEGON2019fractal} the existence of a ``unit ball", and got an estimate that the Minkowski dimension of its fractal boundary is bounded above by $m-\beta$, where $\beta=m(\alpha-(1-\frac{1}{m}))\in(0,1)$. 
\end{remark}

\section{The formation of a plant leaf}

\textcolor{blue}{The idea of using superlevel sets in the study of boundary regularity \cite{xia2011boundary} was partly inspired by the vein structure of tree leaves. We now return to this biological motivation and describe a mathematical model for the dynamic formation of a plant leaf. This example is not merely an application of the previous theory; rather, it was one of the motivations for studying branching transport structures. The purpose of the model is to show how the economy-of-scale principle in ramified optimal transportation can lead naturally to leaf-like shapes and venation patterns.}

In \cite{Xia2007leaf}, a mathematical model was built to understand the 
dynamic formation of a plant leaf. Plant leaves exhibit diverse and 
elaborate shapes and venation patterns. Why do tree leaves grow in such 
an intricate way? What is the mathematics behind it? To answer these 
questions, we first look at the basic functions of leaves. A tree leaf 
transports resources such as water and solutes from the root to its 
tissues via the xylem, absorbs solar energy at its cells through 
photosynthesis (a function of surface area), and then transports the 
chemical products (carbohydrates) synthesized in the leaf back to the root 
through the phloem. Thus, a leaf tends to increase surface area as much as 
possible to maximize metabolic capacity, while also maximizing internal 
efficiency by developing an effective biological transport system.

\textcolor{blue}{From the point of view of optimal transportation, a leaf may be regarded as a transport network connecting one source, the petiole or root node, to a large number of distributed sinks, the cells of the leaf. The model is based on the same mechanism as the ramified transport problem discussed in the previous section. When $0<\alpha<1$, the cost $m^\alpha$ is concave in the transported mass $m$, so transporting mass together along a shared route is cheaper than transporting the same masses separately. This favors branching networks. Thus, the venation pattern of a leaf can be interpreted as an efficient transport system that balances the benefit of producing additional leaf area against the cost of supplying that area with water and nutrients.}

For simplicity, a leaf may be viewed as a finite union of square cells, 
centered on a grid
\[
\Gamma_h = \{ (mh,nh): m,n \in \mathbb{Z} \} \subset \mathbb{R}^2,
\]
with mesh size $h$. A leaf will be a subset of $\Gamma_h$, yet not every 
subset corresponds to a reasonable leaf. One aim is to understand what 
distinguishes realistic leaf structures \textcolor{blue}{from arbitrary collections of grid cells}.

Let
\[
\Omega = \{x_1,x_2,\dots,x_N\} \subset \Gamma_h
\]
be any finite subset representing a prospective leaf. The amount of water 
needed at each cell is proportional to its area ($=h^2$). Without loss of 
generality, we may assume it equals $h^2$. Hence, each $x_i \in \Omega$ 
corresponds to a particle of mass $h^2$ located at $x_i$. The aim is to 
transport these particles to the node $O$ (the point of attachment of the 
leaf to the branch) in a cost-efficient way.

A transport system for the leaf is modeled by a transport path
\[
G = \big( V(G),E(G),w \big)
\]
from $O$ to the measure $\sum_{i=1}^N (h^2)\,\delta_{x_i}$. The cost 
function on such a transport system is a modified version of 
$\mathbf{M}_{\alpha}(G)$. As before, the parameter $\alpha \in [0,1)$ 
\textcolor{blue}{reflects the economy of scale in transportation. Smaller values of $\alpha$ increase the relative advantage of transporting mass along shared routes.} The additional parameter 
$\beta >0$ reflects the principle that when a direction of transport is 
favored (for instance, along a main vein), it is cheaper to move in that 
direction than in others. The cost grows as the deviation from the favored 
direction increases. \textcolor{blue}{Thus $\alpha$ measures the strength of the transport economy of scale, while $\beta$ controls the directional persistence of the veins. Different choices of these two parameters lead to different leaf shapes and different venation patterns.}

\begin{definition}
For $\alpha \in [0,1)$ and $\beta>0$, the cost of a transport path 
$G=(V,E,w)$ is defined by
\[
\mathbf{M}_{\alpha,\beta}(G)
   := \sum_{e \in E(G)} m_\beta(e^+)\, \big(w(e)\big)^{\alpha} 
      \, \text{length}(e),
\]
where the weight factor $m_\beta:V(G)\to (0,\infty]$ is defined 
inductively as follows:
\begin{itemize}
\item $m_\beta(O) = 1$;
\item if $v \in V(G)\setminus \{O\}$ and $p(v)$ is the parent of $v$, then
\[
m_\beta(v) = m_\beta(p(v)) \, H_\beta(\vec e_v,\vec e_{p(v)}),
\]
with
\[
H_\beta(u,v) =
\begin{cases}
|u\cdot v|^{-\beta} = \dfrac{1}{\cos^\beta(\theta)}, & u\cdot v >0, \\
+\infty, & \text{otherwise},
\end{cases}
\]
where $\theta$ is the angle between the unit vectors $u$ and $v$.
\end{itemize}
\end{definition}

Thus $m_\beta$ penalizes sharp turns in transport paths. In biological 
terms, this reflects that flows are cheaper when aligned with existing 
vein directions, but prohibitively costly when reversed. \textcolor{blue}{This directional penalty favors organized venation patterns rather than arbitrary branching networks.}

\textcolor{blue}{The above cost functional keeps the basic feature of ramified transport while incorporating a biological directional effect. The factor $w(e)^\alpha$ represents economy of scale in transportation, while $m_\beta(e^+)$ records the accumulated cost of changing direction along the path from the root to the endpoint. Thus a low-cost network tends to balance efficient branching with directional persistence.}

The growth of a plant leaf is then modeled as a dynamic process of 
generating new cells. At each stage, the leaf develops an optimal transport 
system for existing cells with respect to $\mathbf{M}_{\alpha,\beta}$. Let
\[
\mathcal{A}_h :=
\Big\{ (\Omega,G): \Omega \subset \Gamma_h,\ 
     G \text{ is an optimal transport system of }\Omega
     \text{ under } \mathbf{M}_{\alpha,\beta}\Big\}.
\]

As the environment changes, the leaf may generate new cells at its 
boundary. Each potential new cell produces revenue proportional to its 
area:
\[
\text{revenue} = \epsilon h^2, \qquad \epsilon >0,
\]
while the expense depends on the transport cost of connecting it to the 
existing system. \textcolor{blue}{Thus the model treats leaf growth as a local competition between metabolic gain and transportation expense.}

Given $(\Omega,G) \in \mathcal{A}_h$, let $B$ be the boundary of $\Omega$. 
For $x \in \Gamma_h\setminus\Omega$ and $b\in B$, the transport cost of 
$x$ via $b$ is defined by
\[
C_\Omega(x,b) := h^{2\alpha}\,|x-b|\, m_\beta(b)\,
   H_\beta\!\left(\frac{x-b}{|x-b|},\vec e_b\right)
   + P_G(h^2,b),
\]
where $P_G(y,v)$ measures the incremental cost of adding mass $y$ at 
$v \in V(G)$:
\[
P_G(y,v) = \sum_{u \in P_v \setminus \{O\}}
   m_\beta(u)\Big[ \big(w(e_u)+y\big)^{\alpha}
    - w(e_u)^\alpha\Big] \, \text{length}(e_u).
\]
The expense for generating a new cell at $x$ is then
\[
C_\Omega(x) := \min_{b\in B} C_\Omega(x,b).
\]

\textbf{Selection principle:}
\textit{
A new cell $x$ is generated if and only if
\[
C_\Omega(x) \leq \epsilon h^2.
\]
}

\textcolor{blue}{This selection principle has a direct biological interpretation. A new cell is added precisely when the revenue gained from increasing photosynthetic area is at least as large as the additional transportation cost needed to support that cell. Since the cost $C_\Omega(x)$ depends on the current transport path $G$, the existing venation structure influences future growth. Conversely, once new cells are added, the transport demand changes and the optimal venation pattern must be updated. Thus the model couples shape formation with network optimization.}

This leads to an updated set $\tilde\Omega$ and transport path $\tilde G$, 
and hence to a growth operator
\[
L_{\epsilon,h}:\mathcal{A}_h \to \mathcal{A}_h, \qquad
L_{\epsilon,h}(\Omega,G) = (\tilde\Omega,\tilde G).
\]

Iterating $L_{\epsilon,h}$ eventually stabilizes: for 
$\alpha \in (1/2,1)$, the leaf remains confined to a bounded region. \textcolor{blue}{Indeed, Proposition~3.2 in \cite{Xia2007leaf} gives a radius $R$, independent of $h$, such that the growth operator maps leaves contained in $B_R(O)$ back into $B_R(O)$. Since $B_R(O)\cap\Gamma_h$ is finite, the increasing sequence of generated cell sets must eventually stabilize.} One 
thus defines the \emph{generation map}
\[
g_{\epsilon,h}(\Omega,G) = \lim_{n\to\infty} L_{\epsilon,h}^n(\Omega,G).
\]

\begin{definition}
For $\epsilon>0$ and $h>0$, a pair $(\Omega,G)\in\mathcal{A}_{\epsilon,h}$ 
is called an \emph{$(\epsilon,h)$-leaf} if it can be obtained from the bud
\[
\Omega_0=\{O\}, \quad G_0=(\{O\},\emptyset,-),
\]
by successive application of generation maps with strictly increasing 
revenue parameters $0<\epsilon_1<\cdots<\epsilon_k=\epsilon$.
\end{definition}

The preceding definitions lead to the following schematic numerical procedure for generating a leaf.

\begin{algorithm}[H]
\caption{Leaf Growth--Relaxation Algorithm}
\label{alg:leaf_growth}

\KwIn{
Mesh size $h$, parameters $\alpha\in[0,1)$ and $\beta>0$, and revenue levels
$0<\epsilon_1<\epsilon_2<\cdots<\epsilon_k=\epsilon$
}

\KwOut{
An $(\epsilon,h)$-leaf $(\Omega,G)$
}

\BlankLine

\textbf{Initialization:}\\
Set
\[
\Omega\leftarrow \{O\},
\qquad
G\leftarrow (\{O\},\emptyset,-).
\]

\BlankLine

\For{$j=1,\ldots,k$}{

\Repeat{no new cells are added}{

\textbf{Candidate cell selection:}\\
Let $B$ be the boundary of $\Omega$.\\
Define
\[
D(\Omega)
=
\{x\in\Gamma_h\setminus\Omega:
x \text{ is adjacent to } \Omega\}.
\]

\BlankLine

\textbf{Expense evaluation:}\\
\ForEach{$x\in D(\Omega)$}{
Compute
\[
C_\Omega(x)
=
\min_{b\in B} C_\Omega(x,b).
\]
}

\BlankLine

\textbf{Growth step:}\\
Select
\[
S
=
\{x\in D(\Omega):
C_\Omega(x)\leq \epsilon_j h^2\}.
\]

\If{$S=\emptyset$}{
Stop the iteration at revenue level $\epsilon_j$.
}

Set
\[
\Omega\leftarrow \Omega\cup S.
\]

\BlankLine

\textbf{Initial graph update:}\\
Connect each newly selected cell $x\in S$ to the existing transport network through a boundary point
\[
b_x\in\arg\min_{b\in B} C_\Omega(x,b).
\]
Let $G$ denote this initially updated graph.

\BlankLine

\textbf{Local transport relaxation:}\\
\Repeat{no admissible parent change decreases $\mathbf M_{\alpha,\beta}(G)$}{
\ForEach{$u\in V(G)\setminus\{O\}$}{
Search over admissible new parents $v$ for $u$.\\
Let $G_{u\to v}$ be the graph obtained by changing the parent of $u$ to $v$.\\
\If{$\mathbf M_{\alpha,\beta}(G_{u\to v})<\mathbf M_{\alpha,\beta}(G)$}{
Set $G\leftarrow G_{u\to v}$.
}
}
}

}

}

\Return{$(\Omega,G)$}

\end{algorithm}

\textcolor{blue}{Algorithm~\ref{alg:leaf_growth} gives a schematic description of the numerical procedure used to generate a leaf. Starting from the bud, the algorithm increases the environmental revenue parameter and alternates between cell selection and local relaxation of the transport network. The output is an $(\epsilon,h)$-leaf in the sense of the definition above. The relaxation step is energy-decreasing, but the procedure should be viewed as a local numerical approximation rather than a certified global minimization at every stage. If $N=|V(G)|$ and $M=|D(\Omega)|$, evaluating all candidate cells costs $O(M|B|)$ in a direct implementation, while a naive parent-reassignment sweep costs $O(N^2)$. With $r$ relaxation sweeps, one growth step has naive cost $O(M|B|+rN^2)$. Restricting the parent search to $k_0$ nearby admissible vertices reduces the practical relaxation cost to about $O(rk_0N)$.}

\begin{figure}[h]
\centering
\includegraphics[width=5cm]{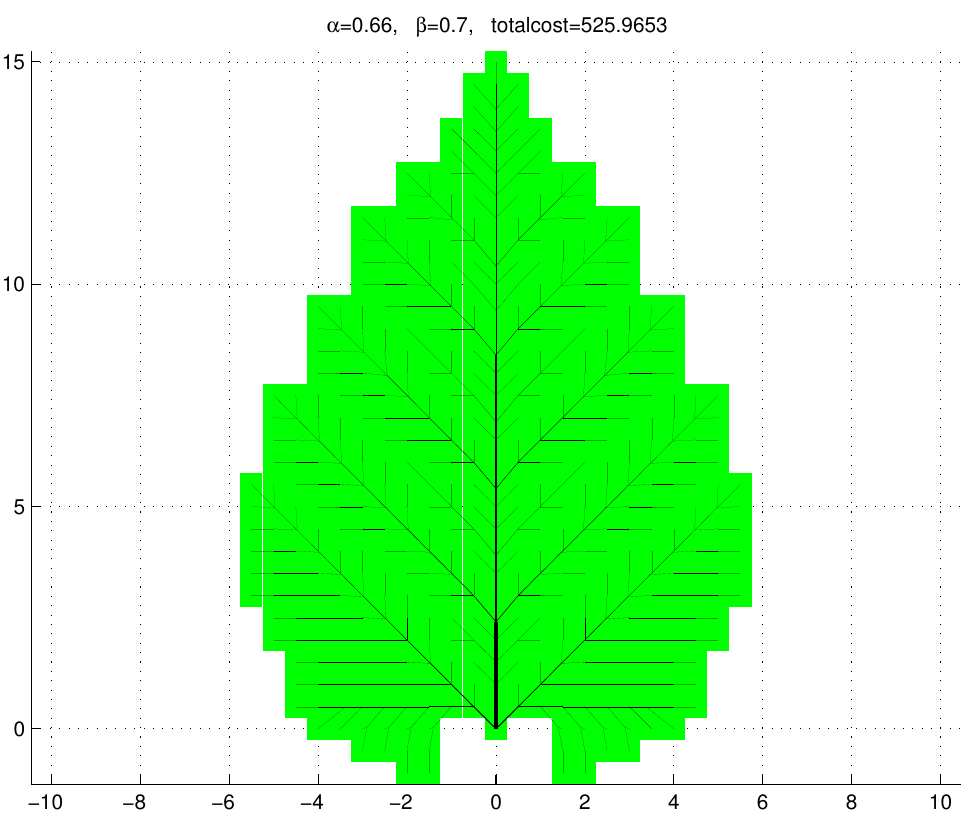}\quad
\includegraphics[width=5cm]{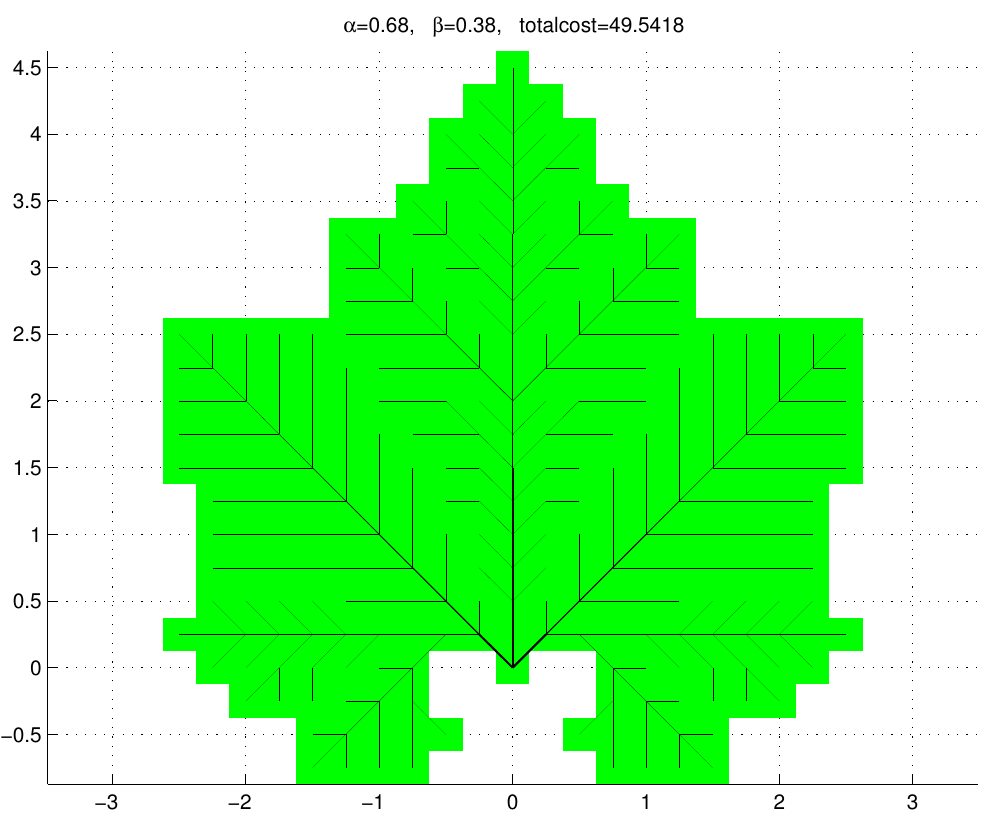}
\caption{\textcolor{blue}{Numerically generated leaves under the dynamic ramified transport growth model. The simulations illustrate how transport economy, directional persistence, and environmental revenue interact to produce leaf-like shapes and branching venation patterns, including examples resembling maple and mulberry leaves.}}
\label{fig:leaves}
\end{figure}

As a result, the final shape and venation pattern of a leaf are determined 
by the transport cost function and the environment. \textcolor{blue}{More precisely, they are produced through a dynamic growth process in which new cells are selected and the transport network is repeatedly updated.} Numerical simulations 
based on this model produce leaves resembling those of maple and mulberry 
(Figure~\ref{fig:leaves}). \textcolor{blue}{The branching veins are not prescribed in advance; rather, they emerge from the repeated optimization of the transport system as new cells are added.}

\textcolor{blue}{The model also has a limiting interpretation as the mesh size $h$ tends to zero. In the limit, the discrete leaf may be represented by a Radon measure $\mu$, while its transport system is represented by a vector measure $\Theta$ satisfying}
\[
\operatorname{div}\Theta
=
M(\mu)\delta_O-\mu.
\]
\textcolor{blue}{This connects the discrete growth model with the measure-theoretic framework used elsewhere in this article and provides another link between the numerical simulations and the geometric-measure-theoretic analysis of branching transport networks.}

\section{A unified theory of landscape functions}

\textcolor{blue}{After studying properties of the $\mathbf{M}_\alpha$ mass, it is natural to consider its first variation, or equivalently, the marginal cost of transporting an additional infinitesimal amount of mass along a given network. This point of view leads to the notion of a landscape function. Roughly speaking, a landscape function records the accumulated marginal transportation cost from a source to points of the transport network.}

\textcolor{blue}{The idea of using marginal transportation cost in ramified optimal transportation appears, for instance, in applications such as leaf formation \cite{Xia2007leaf} and diffusion-limited aggregation driven by optimal transportation \cite{xia2007DLA}. A systematic study was initiated by Santambrogio \cite{Santa}, who introduced landscape functions through a Lagrangian formulation of irrigation plans. Santambrogio's landscape function is associated with an optimal irrigation plan starting from a single source. Later, Xia introduced an Eulerian version for finite acyclic transport paths between atomic measures \cite{xia2014landscape}. This formulation allows multiple sources and targets, and the transport path need not be optimal. In Schiffman's dissertation \cite{Schiffman2021}, these two approaches were unified using geometric measure theory by defining landscape functions for a broader class of $\alpha$-balanced rectifiable currents.}

\textcolor{blue}{This section reviews these three related approaches. We first describe Santambrogio's Lagrangian formulation for optimal irrigation plans from a single source \cite{Santa}. We then describe Xia's Eulerian formulation for finite acyclic transport paths between atomic measures, where the landscape equation can be written using the incidence matrix \cite{bapat} of the underlying graph \cite{xia2014landscape}. Finally, we explain Schiffman's geometric-measure-theoretic framework, in which both previous notions arise as special cases of landscape functions for $\alpha$-balanced rectifiable currents \cite{Schiffman2021}. Thus the purpose of this section is to show how landscape functions connect first variation, graph-theoretic potentials, and the current-theoretic formulation of ramified transport.}

\subsection{Santambrogio's approach of landscape function}\label{sec: Santambrogio}

We begin with Santambrogio's Lagrangian approach. In the study of river basins, geophysics yields problems very similar to ramified optimal transportation.
While studying the
configuration of a river basin, the main objects are the landscape
function giving the altitude of any point of the region and a river
channel network which is the datum of all the streams that concur to
bring water to lakes. In a seminal paper \cite{Santa} of Santambrogio, he
introduced the concept ``landscape function" associated with an
optimal transport path from a Dirac mass (representing a single
lake) to another probability measure (representing the region).

His approach was described in a Lagrangian way. Here, we quickly set the Lagrangian framework of branched transport and its main features. For more details, we refer to the book \cite{bernot2008optimal}.

Let us start with the description of the irrigation problem. We denote by $\Gamma(\mathbb{R}^d)$ the set of $1$-Lipschitz curves in $\mathbb{R}^d$ parameterized on $[0,\infty]$, endowed with the topology of uniform convergence on compact sets.

We call \emph{irrigation plan} any probability measure $\eta \in Prob({\Gamma})$ satisfying the following finite-length condition
\begin{equation}\label{fin_len}
\mathbf{L}(\eta) :=\int_\Gamma L(\gamma) \:d \eta(\gamma) < +\infty,
\end{equation}
where $L(\gamma) = \int_0^\infty \lvert{\dot{\gamma}(t)}\rvert d t$. Notice that any irrigation plan is concentrated on $\Gamma^1(\mathbb{R}^d) :=\{ \gamma : L(\gamma) < \infty\}$. We denote by $\mathrm{IP}(\mathbb{R}^d)$ the set of all irrigation plans $\eta\in Prob({\Gamma})$. If $\mu$ and $\nu$ are two probability measures on $\mathbb{R}^d$, one says that $\eta \in \mathrm{IP}(\mathbb{R}^d)$ irrigates $\nu$ from $\mu$ if one recovers the measures $\mu$ and $\nu$ by sending the mass of each curve respectively to its initial point and to its final point, which means that
\[
(\pi_0)_\#\eta = \mu \text{ and }
(\pi_\infty)_\#\eta = \nu,
\]
where $\pi_0(\gamma) = \gamma(0)$, $\pi_\infty(\gamma) = \gamma(\infty) :=\lim_{t \to +\infty} \gamma(t)$ and $f_\#\eta$ denotes the push-forward of $\eta$ by $f$ whenever $f$ is a Borel map. We denote by $\mathrm{IP}(\mu,\nu)$ the set of irrigation plans irrigating $\nu$ from $\mu$:
\[\mathrm{IP}(\mu,\nu) = \{ \eta \in \mathrm{IP}(\mathbb{R}^d) : (\pi_0)_\#\eta = \mu, (\pi_\infty)_\#\eta = \nu\}.\]
If $\eta$ is a given irrigation plan, we define the multiplicity at $x$, that is the total mass passing by $x$, as
\[\theta_\eta(x) = \eta(\{\gamma \in \Gamma : x \in \gamma\}),\]
where $x\in \gamma$ means that $x$ belongs to the image of the curve $\gamma$. Finally, for any nonnegative function $f$, we denote by $\int_\gamma f(x) d x$ the line integral of $f$ along $\gamma \in \Gamma$:
\[\int_\gamma f(x) d x :=\int_0^{+\infty} f(\gamma(t)) \lvert{\dot{\gamma}(t)}\rvert d t.\]

For $\alpha \in [0,1]$, we consider the \emph{irrigation cost} $\mathbf{I}_\alpha : \mathrm{IP}(\mathbb{R}^d) \to [0,\infty]$ defined by
\[\mathbf{I}_\alpha(\eta) :=\int_\Gamma \int_\gamma \theta_\eta(x)^{\alpha -1} d xd\eta(\gamma),\]
with the conventions $0^{\alpha -1} = \infty$ if $\alpha <1$, $0^{\alpha -1} = 1$ otherwise, and $\infty \times 0 = 0$. If $\mu,\nu$ are two probability measures on $\mathbb{R}^d$, the irrigation (or branched transport) problem consists in minimizing the cost $\mathbf{I}_\alpha$ on the set of irrigation plans which send $\mu$ to $\nu$, which reads
\begin{equation}\label{lag_pb}\tag{$\text{LI}_\alpha$}
\min_{\eta \in \mathrm{IP}(\mu,\nu)}\quad \int_\Gamma \int_\gamma \theta_\eta(x)^{\alpha -1} d x d\eta(\gamma).
\end{equation}
We set $Z_\eta(\gamma) = \int_\gamma \theta_\eta(x)^{\alpha-1} d x$ so that the cost may expressed as
\[\mathbf{I}_\alpha(\eta) = \int_\Gamma Z_\eta(\gamma) d\eta(\gamma).\]

We refer to the book \cite{bernot2008optimal} for the equivalence of minimizing $\mathbf{I}_\alpha$ over irrigation plans with minimizing $\mathbf{M}_\alpha$ over transport paths.

The following results are extracted from \cite{xia2003, BCM2005}.

\begin{proposition}[First variation inequality for $\mathbf{I}_\alpha$]
If $\eta$ is an irrigation plan with $\mathbf{I}_\alpha(\eta)$ finite, then for all irrigation plan $\tilde{\eta}$ the following holds:
\begin{equation}
\label{1st_var}
\mathbf{I}_\alpha(\tilde{\eta}) \leq \mathbf{I}_\alpha(\eta) + \alpha \int Z_\eta(\gamma) d(\tilde{\eta}-\eta).
\end{equation}
Notice that the integral $\int Z_\eta d(\tilde{\eta}-\eta)$ is well-defined since $\int Z_\eta d\eta < \infty$ and $Z_\eta$ is nonnegative, though it may be infinite.
\end{proposition}

 We now recall the basic definitions and properties about landscape functions. Given an optimal irrigation plan $\eta \in \mathrm{IP}(\delta_0,\nu)$, we say that a curve $\gamma$ is $\eta$-good if
\begin{itemize}
\item the quantity $Z_\eta(\gamma) = \int_\gamma \theta_\eta(x)^{\alpha-1}  d x$ is finite,
\item for all $t < T(\gamma)$,
\[\theta_\eta(\gamma(t)) = \eta(\{\tilde{\gamma}\in \Gamma(\mathbb{R}^d) : \gamma = \tilde{\gamma} \text{ on } [0,t]\}),\]
where $T(\gamma)=\inf \{t\in [0,\infty]: \gamma(s)=\gamma(\infty) \text{ for all } s\in[ t, \infty]\}$ is the stopping time of $\gamma$.
\end{itemize}

One may prove by optimality that $\eta$ is concentrated on the set of $\eta$-good curves. Moreover it is proven in \cite{Santa} that for all $\eta$-good curves $\gamma$, the quantity $Z_\eta(\gamma)$ depends only on the final point $\gamma(\infty)$ of the curve, thus we may define the landscape function $z_\eta$ as follows:

\begin{definition}
When $\eta \in IP(\delta_O,\nu)$ is $\alpha$-optimal,
the landscape function $z_\eta : \mathbb{R}^d \rightarrow [0,+\infty]$ is defined by
	\[
z_\eta(x)
		:= \left \{ \begin{array}{ll}
		Z_\eta(\gamma), & \text{ if }\gamma \text{ is an } \eta\text{-good curve s.t. }x = \gamma(\infty),\\
		+\infty, & \text{if no $\eta$-good curve ends at } x	\end{array}\right.
	\]
\end{definition}

Notice that for an optimal $\eta$ the cost may be expressed in terms of $z_\eta$:
\[\mathbf{I}_\alpha(\eta) = \int_\Gamma Z_\eta(\gamma) d\eta(\gamma) = \int_{\mathbb{R}^d} z_\eta(x) d\nu(x).\]
Finally, one may show that $z_\eta$ is lower semicontinuous and that the inequality $z_\eta(x)\geq |x|$ holds.

\subsection{Xia's approach of landscape function}
Santamborgio's approach works on optimal irrigation plan from a single source. 
Let $\mathbf{a}=\delta _{S}$ be the initial source and the target $\mathbf{b}$ be an atomic measure on $X$.
Let $P=\left\{ V\left( P\right) ,E\left( P\right) ,w\right\}$ be an $\alpha$-optimal transport
path from $\mathbf{a}$ to $\mathbf{b}$. Then for any $x$ in the support of $P
$, there is a unique polyhedral curve $\gamma _{x}$ on $P$ from the initial
source $S$ to $x$. In this case, Santambrogio's landscape function $z\left(
x\right) $ associated with $P$ can be simply written as 
\begin{equation*}
z\left( x\right) =\int_{\gamma _{x}}\theta \left( s\right) ^{\alpha -1}d%
\mathcal{H}^{1}\left( s\right)
\end{equation*}%
where $\mathcal{H}^{1}$ represents the $1-$dimensional Hausdorff measure,
and $\theta(s)$ represents the mass flowing through the point $\gamma _{x}(s)$.

Intuitively,  it is the marginal transportation cost of the mass from the  source to the point. Independently, I 
have also used analogous ideas of the landscape function in
applications such as modeling diffusion-limited
aggregation driven by optimal transportation \cite{xia2007DLA} and analyzing optimal assignment maps \cite{xia2013allocation}.

When consider marginal transport cost, there is no need to restrict on single source and optimal transport path. So, I would like to generalize it using geometric measure theory.
I was stuck with workable approach for general case, so I provided a version for discrete case in the article
\cite{xia2014landscape}, where we generalized the concept of landscape function by
allowing multiple sources rather than a single common source, and the transport path is not necessarily optimal.

We now describe this Eulerian approach as follows.

\begin{definition}
Given two atomic measures $\mathbf{a,b}$ of equal mass, let $$P=\left\{ V\left( P\right) ,E\left( P\right) ,\mathbf{\omega }\right\} $$
be an acyclic transport path from $\mathbf{a}$ to $\mathbf{b}$ with $V(P)=\{v_1,v_2,\cdots,v_J\}$, $E(P)=\{e_1,e_2,\cdots, e_K\}$, 
and $0\leq
\alpha <1$. A function $Z:V\left( P\right) \rightarrow \mathbb{R}$ is called
an $\alpha -$landscape function associated with $P$ if for each edge $e\in
E\left( P\right) $, it holds that 
\begin{equation}
Z\left( e^{+}\right) -Z\left( e^{-}\right) =\mathbf{\omega }\left( e\right)
^{\alpha -1}length\left( e\right) \text{.}  \label{landscape_eqn}
\end{equation}
\end{definition}

To solve (\ref{landscape_eqn}), we first express it in matrix notation. Let $Q_P $ be the
incidence matrix of the underlying graph $G_P$, 
\[%
\boldsymbol{w}\mathbf{=}\left[ \mathbf{\omega }\left( e_{1}\right) ,\cdots ,%
\mathbf{\omega }\left( e_{K}\right) \right] \text{ and }\boldsymbol{c}%
\mathbf{=}\left[ c\left( v_{1}\right) ,\cdots ,c\left( v_{J}\right) \right] \]
with 
\begin{equation*}
c\left( v\right) =\left\{ 
\begin{array}{c}
-m_{i},\text{\ if }v=x_{i}\text{\ for some }i=1,...,k \\ 
n_{j},\text{\ if }v=y_{j}\text{\ for some }j=1,...,\ell \\ 
0,\text{\ otherwise }%
\end{array}%
\right. \end{equation*}%
for each $v\in V\left( P\right) =\left\{ v_{1},v_{2},\cdots
,v_{J}\right\} $.
  Then, in matrix notation, the balance equation (\ref{balance}) can be expressed
as 
\begin{equation}
\boldsymbol{w}Q_P ^{\prime }\mathbf{=}\boldsymbol{c}
\label{path_matrix}
\end{equation}%
where $Q_P ^{\prime }$ denotes the transpose of $Q_P$. 
When the underlying graph of  $%
P$ is acyclic and connected, then $\boldsymbol{w}$ is uniquely determined
by the incidence matrix $Q_P$ and $\boldsymbol{c}$.

Let $Z=\left[ Z\left( v_{1}\right) ,Z\left( v_{2}\right) ,\cdots ,Z\left(
v_{J}\right) \right] $ and 
\[
D=[\mathbf{\omega }\left( e_{1}\right) ^{\alpha -1}length\left( e_{1}\right)
,\cdots ,\mathbf{\omega }\left( e_{K}\right) ^{\alpha -1}length\left(
e_{K}\right) ].
\] Then, in matrix notations, (\ref{landscape_eqn}) becomes
\begin{equation}
ZQ_P =D .\label{landscape_system}
\end{equation}%

From linear algebra, one can show that 
\begin{itemize}
\item the dimension of the solution space of the system
(\ref{landscape_system}) is equal to the connected components of the underlying graph of $P$.
\item  In particular, for a connected
transport path $P$, its $\alpha $-landscape function is unique up to a
constant. This agrees with the motivation that the landscape function
represents the elevation of the landscape.
\end{itemize}
\begin{example}
Consider the transport path $P$ illustrated in Figure \ref{landscape_figure}.
\begin{figure}[h]
\centering
\includegraphics[width=0.5\textwidth]{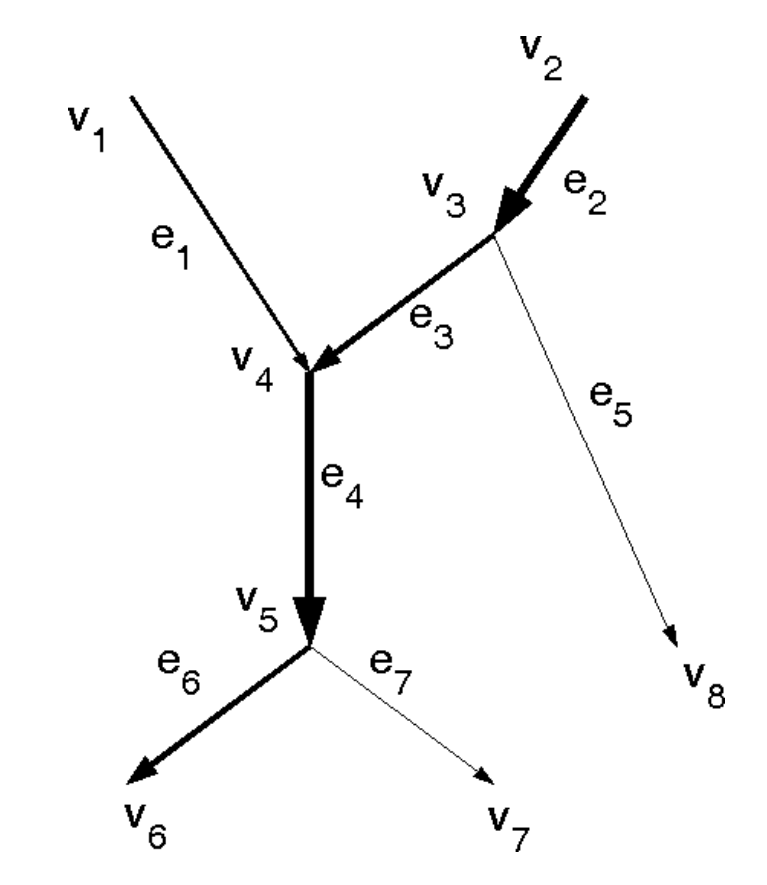}
\caption{Landscape functions associated with a transport path.}
\label{landscape_figure}
\end{figure}
Then, the incidence matrix of $G_{P}$ is 
\[
Q\left( G_{P}\right) =
\begin{bmatrix}
-1 & 0 & 0 & 0 & 0 & 0 & 0 \\ 
0 & -1 & 0 & 0 & 0 & 0 & 0 \\ 
0 & 1 & -1 & 0 & -1 & 0 & 0 \\ 
1 & 0 & 1 & -1 & 0 & 0 & 0 \\ 
0 & 0 & 0 & 1 & 0 & -1 & -1 \\ 
0 & 0 & 0 & 0 & 0 & 1 & 0 \\ 
0 & 0 & 0 & 0 & 0 & 0 & 1 \\ 
0 & 0 & 0 & 0 & 1 & 0 & 0%
\end{bmatrix}
\]

Let 
\begin{eqnarray*}
\mathbf{a} &=&\frac{1}{4}\delta _{v_{1}}+\frac{3}{4}\delta _{v_{2}} \\
\mathbf{b} &=&\frac{1}{2}\delta _{v_{6}}+\frac{1}{4}\delta _{v_{7}}+\frac{1}{%
4}\delta _{v_{8}}
\end{eqnarray*}%
The weight function $\boldsymbol{w}$ solving (\ref%
{path_matrix}) is given by 
\begin{equation*}
\boldsymbol{w}=[1/4,3/4,1/2,3/4,1/4,1/2,1/4].
\end{equation*}%

Assume the coordinates of vertices $\left\{
v_{i}\right\} $ in $\mathbb{R}^{2}$ are given by the following table:%
\begin{equation*}
\begin{tabular}{|l|l|l|l|l|l|l|l|l|}
\hline
vertex $v=\left( x,y\right) $ & $v_{1}$ & $v_{2}$ & $v_{3}$ & $v_{4}$ & $%
v_{5}$ & $v_{6}$ & $v_{7}$ & $v_{8}$ \\ \hline
$x$ & $3$ & $5.5$ & $5$ & $4$ & $4$ & $3$ & $5$ & $6$ \\ \hline
$y$ & $9$ & $9$ & $8$ & $7$ & $5$ & $4$ & $4$ & $5$ \\ \hline
\end{tabular}%
\end{equation*}

Then, the lengths of edges are given by the vector 
\begin{equation*}
L=[\sqrt{5},\sqrt{5}/2,\sqrt{2},2,\sqrt{10},\sqrt{2},\sqrt{2}].
\end{equation*}%

Now, we can find landscape functions associated with $P$ by solving the
system (\ref{landscape_system}). For instance, 

\begin{itemize}
\item when $\alpha =0.5$, the corresponding $\alpha $-landscape function is
\begin{equation*}
Z=[0,1.1811,2.4721,4.4721,6.7815,8.7815,9.6100,8.7967]+c\vec{1}
\end{equation*}%
\item When $%
\alpha =0.85$, the corresponding $\alpha $-landscape function is
\begin{equation*}
Z=[0,0.0164,1.1838,2.7529,4.8411,6.4103,6.5822,5.0770]+c\vec{1}
\end{equation*}%
\end{itemize}
These landscape functions are unique up to a constant.
\end{example}

Using the values of $Z$ on vertices, one may linearly extend it to points on edges, leads to a (generalized) landscape function $Z: spt(P)\rightarrow \mathbb{R}$. 

We now study properties of landscape functions. First, landscape functions give another formula for the transport cost $\mathbf{M}_{\alpha}(P)$:
\begin{proposition}
For any $\alpha -$landscape function $Z$ associated with a transport path $%
P\in Path\left( \mathbf{a,b}\right) $, it holds that 
\begin{equation*}
\int_{X}Zd\left( \mathbf{b-a}\right) =\mathbf{M}_{\alpha }\left( P\right) 
\text{.}
\end{equation*}%
In particular, if $P$ is an $\alpha$-optimal transport path from $\mathbf{a}$
to $\mathbf{b}$, then 
\begin{equation*}
\int_{X}Zd\left( \mathbf{b-a}\right) =d_{\alpha }\left( \mathbf{a,b}\right) 
\text{.}
\end{equation*}
\end{proposition}

It turns out that the study of $\alpha$-landscape functions is closely related to the study of $p$-harmonic functions on directed graphs where $p=\frac{\alpha}{\alpha-1}$ is \textit{non-positive}. More precisely, let $G$ be an acyclic directed graph with a vertex set $V\left( G\right) $
and an edge set $E\left( G\right) $ of directed edges.

\begin{definition}
For any two vertices $v,\tilde{v}\in V\left( G\right) $, we say $v\prec 
\tilde{v}$ if there exists a list of edges $\left\{
e_{i_{1},}e_{i_{2}},\cdots ,e_{i_{k}}\right\} $ in $E\left( G\right) $ with $%
e_{i_{h}}^{+}=e_{i_{h+1}}^{-}$ for $h=1,2,\cdots ,k-1$, and $e_{i_{1}}^{-}=v$%
, $e_{i_{k}}^{+}=\tilde{v}$.
\end{definition}

Note that since $G$ is acyclic, the partial order $\prec $ is well-defined
on $V\left( G\right) $. In particular, $e^{-}$ $\prec $ $e^{+}$ for any edge 
$e\in E\left( G\right) $.

\begin{definition}
\label{mototone_def}Let $\tilde{V}$ be a subset of $V\left( G\right) $. A
function $u:\tilde{V}\rightarrow \mathbb{R}$ is monotone increasing with
respect to $G$ if for any $x,y\in \tilde{V}$ with $x\prec $ $y$, it holds
that $u\left( x\right) <u\left( y\right) $.
\end{definition}

Let $\mathcal{F}_{G}$ be the family of all monotone increasing functions $%
u:V\left( G\right) \rightarrow \mathbb{R}$ with respect to $G$. For
instance, any landscape function $Z$ associated with an acyclic transport
path $P$ is monotone increasing with respect to the underlying graph $G_{P}$.

For each $u\in \mathcal{F}_{G}$ and $p\leq 0$, the $p$\textit{-energy} of $u$
is given by 
\begin{equation}
\mathbf{E}_{p}\left( u\right) =\sum_{e\in E\left( G\right) }\left\vert
\nabla u\left( e\right) \right\vert ^{p}length\left( e\right)
\label{p-energy}
\end{equation}%
where for each $e\in E\left( G\right) $ with unit directional vector $\vec{e}
$, 
\begin{equation*}
\nabla u\left( e\right) =\frac{u\left( e^{+}\right) -u\left( e^{-}\right) }{%
length\left( e\right) }\vec{e}.
\end{equation*}

We are interested in minimizing $\mathbf{E}_{p}\left( u\right)$ among $u\in 
\mathcal{F}_{G}$ with a Dirichlet boundary condition. Let $\partial G$ be a subset of $V\left( G\right) $ such that 
\begin{equation*}
\left\{ 
\begin{array}{c}
v\in V\left( G\right) :\text{ either there is no edge }e\in E\left( G\right) 
\text{ with }e^{+}=v \\ 
\text{or there is no edge }e\in E\left( G\right) \text{ with }e^{-}=v%
\end{array}%
\right\} \subseteq \partial G\subseteq V\left( G\right) .
\end{equation*}%
In other words, $\partial G$ contains all source and sink vertices in $%
V\left( G\right) $, and may contain some other vertices. We view $\partial G$
as the {\it boundary set }of $V\left( G\right) $.

We now consider the following Dirichlet problem on graphs.\\
\begin{problem}
Suppose $u_{0}:\partial G\rightarrow \mathbb{R}$ is monotone increasing with
respect to $G$, and $p<0$. Minimize 
\begin{equation*}
\mathbf{E}_{p}\left( u\right) =\sum_{e\in E\left( G\right) }\left\vert
\nabla u\left( e\right) \right\vert ^{p}length\left( e\right)
\end{equation*}%
among all $u\in U=\left\{ u:V\left( G\right) \rightarrow \mathbb{R}|u\in 
\mathcal{F}_{G}\text{ and }u=u_{0}\text{ on }\partial G\right\} .$
\end{problem}

By identifying $u\in\mathcal{F}_{G}$ with the value $(u(v_1), u(v_2), \cdots, u(v_I))$, it turns out that this minimization problem is equivalent to a minimization of some strictly convex function $f$ on a nonempty convex domain $\Omega$ in $\mathbb{R}^I$.  The ``nonempty convex" property of $\Omega$ is guaranteed by the monotonicity of $u_0$. The associated Euler-Langrange equation is 

\begin{equation}
\sum_{\substack{ e\in E\left( G\right)  \\ ^{e^{+}=v_{i}}}}\left\vert \nabla
u\left( e\right) \right\vert ^{p-1}=\sum_{\substack{ e\in E\left( G\right) 
\\ e^{-}=v_{i}}}\left\vert \nabla u\left( e\right) \right\vert ^{p-1}
\label{vertex_balance_for_minimizer}
\end{equation}%
at each $v_{i}\in V\left( G\right) \setminus \partial G$. Using the discrete
version of the divergence notation: 
\begin{equation*}
\text{div}\left( \vec{V}\right) \left( v\right) =\sum_{\substack{ \text{%
either }e^{+}=v  \\ \text{or }e^{-}=v}}\vec{V}\left( e\right) \cdot \vec{e}
\end{equation*}%
for any vector field $\vec{V}:E\left( G\right) \rightarrow \mathbb{R}^{m}$,
equation (\ref{vertex_balance_for_minimizer}) can be expressed as the $p-$%
Laplace equation on the graph $G:$%
\begin{equation}
\text{div}\left( \left\vert \nabla u\right\vert ^{p-2}\nabla u\right) \left(
v\right) =0  \label{p_laplacian_eqn}
\end{equation}%
for any $v\in V\left( G\right) \backslash \partial G$.

A solution to the $p-$Laplace equation (\ref{p_laplacian_eqn}) in $U$ is
called a $p-$\textit{harmonic function} on the graph $G$. By the strict
convexity of $f$, any $p-$harmonic function $u$ is an $\mathbf{E}_{p}$-
minimizer in $\mathcal{F}_{G}$ with respect to its boundary datum.

The following theorem describes an equivalent relationship between landscape functions and $p$-harmonic maps:
\begin{theorem}
Let $\alpha \in \lbrack 0,1)$ and $p=\frac{\alpha }{\alpha -1}\leq 0$ be the
conjugate of $\alpha $. 
\begin{itemize}
\item[1.] If $Z$ is an $\alpha $-landscape function associated
with an acyclic transport path $P$, then $Z$ is a $p$-harmonic function on
the underlying graph $G_{P}$. Moreover, $\mathbf{M}_{\alpha }\left( P\right)
=\mathbf{E}_{p}\left( Z\right) $.
\item[2.] Conversely, if $u$ is a $p-$harmonic function on
an acyclic graph $G$, then $u$ is an $\alpha $-landscape function associated
with an acyclic transport path $P$ with $%
G$ being its underlying graph. Moreover, $\mathbf{E%
}_{p}\left( u\right) =\mathbf{M}_{\alpha }\left( P\right) $.
\end{itemize}

\end{theorem}

\begin{proof}[Proof sketch]
\textcolor{blue}{The following argument sketches the equivalence proved in \cite{xia2014landscape}. Suppose first that $Z$ is an $\alpha$-landscape function associated with an acyclic transport path $P$. By definition, for each edge $e\in E(P)$,}
\[
Z(e^+)-Z(e^-)
=
\omega(e)^{\alpha-1}\operatorname{length}(e).
\]
\textcolor{blue}{Therefore,}
\[
|\nabla Z(e)|=\omega(e)^{\alpha-1}.
\]
\textcolor{blue}{Since}
\[
p=\frac{\alpha}{\alpha-1},
\qquad
p-1=\frac{1}{\alpha-1},
\]
\textcolor{blue}{we obtain}
\[
|\nabla Z(e)|^{p-1}=\omega(e).
\]
\textcolor{blue}{Thus the $p$-harmonic vertex equation}
\[
\sum_{e^+=v}|\nabla Z(e)|^{p-1}
=
\sum_{e^-=v}|\nabla Z(e)|^{p-1}
\]
\textcolor{blue}{is exactly the mass-balance equation for the transport path $P$ at every interior vertex $v$. Hence $Z$ is $p$-harmonic on the underlying graph $G_P$. Moreover,}
\[
\mathbf E_p(Z)
=
\sum_{e\in E(P)}|\nabla Z(e)|^p\operatorname{length}(e)
=
\sum_{e\in E(P)}\omega(e)^\alpha\operatorname{length}(e)
=
\mathbf M_\alpha(P).
\]

\textcolor{blue}{Conversely, suppose that $u$ is a $p$-harmonic function on an acyclic directed graph $G$. Define a weight on each edge by}
\[
\omega(e)=|\nabla u(e)|^{p-1}.
\]
\textcolor{blue}{The $p$-harmonic equation gives the balance condition for these weights at each interior vertex. Hence the weighted graph defines a transport path $P$ between two atomic measures determined by the net imbalance at the boundary vertices. Since}
\[
\omega(e)^{\alpha-1}
=
|\nabla u(e)|,
\]
\textcolor{blue}{we have}
\[
u(e^+)-u(e^-)
=
\omega(e)^{\alpha-1}\operatorname{length}(e).
\]
\textcolor{blue}{Thus $u$ is an $\alpha$-landscape function associated with the resulting transport path $P$. The same calculation gives}
\[
\mathbf E_p(u)=\mathbf M_\alpha(P).
\]
\end{proof}

We now state some properties of landscape functions that associated with optimal transport paths.
\begin{proposition}
Suppose $P\in Path\left( \mathbf{a,b}\right) $ is an $\alpha $-optimal
transport path. For any $\alpha $-landscape function $Z$ associated with $P,$
it holds that 
\begin{equation}
\sum_{e^{+}=v}\mathbf{\omega }\left( e\right) \nabla Z\left( e\right)
=\sum_{e^{-}=v}\mathbf{\omega }\left( e\right) \nabla Z\left( e\right)
\label{optimal_eqn_1}
\end{equation}%
at each vertex $v\in V\left( P\right) \setminus \left\{
x_{1},\cdots ,x_{k},y_{1},\cdots ,y_{\ell }\right\} $, and 
\begin{equation}
\sum_{e^{+}=v}\left\vert \nabla Z\left( e\right) \right\vert ^{p-1}\nabla
Z\left( e\right) =\sum_{e^{-}=v}\left\vert \nabla Z\left( e\right)
\right\vert ^{p-1}\nabla Z\left( e\right)  \label{optimal_eqn_2}
\end{equation}%
for $p=\frac{\alpha }{\alpha -1}$.
\end{proposition}
\begin{remark}
\label{momentum_remark} For any transport path $P$, by (\ref{balance}),
mass is conserved at each interior vertex of $P$. When the transport path $P$
is optimal, by (\ref{optimal_eqn_1}), momentum is also conserved at each
interior vertex when viewing $\nabla Z(e)=w(e)^{\alpha-1}\vec{e}$ as the
velocity of a moving particle of mass $w(e)$ on each $e$.
\end{remark}
The following proposition says that any $\alpha $-landscape function
associated with an $\alpha -$optimal transport path is Lipschitz. Before
stating the proposition, we first extend the domain of a landscape function $%
Z:V\left( P\right) \rightarrow \mathbb{R}$ to the support of $P$ by linear
extensions of $Z$ on each edge of $P$.
\begin{proposition}
\label{lipschitzcont} Suppose $P\in Path\left( \mathbf{a,b}\right) $ is an $%
\alpha $-optimal transport path for some $\alpha \in \left( 0,1\right) $.
Let $Z$ be an $\alpha $-landscape function associated with $P$. Then, for
any $x,y$ on the same connected components of the support of $G_{P}$, it
holds that 
\begin{equation*}
\left\vert Z\left( x\right) -Z\left( y\right) \right\vert \leq \frac{1}{%
\alpha }\sigma ^{\alpha -1}\left\Vert x-y\right\Vert
\end{equation*}%
where $\sigma =\min_{e\subseteq \gamma }\mathbf{\omega }\left( e\right) $,
and $\gamma $ is the unique curve on $P$ from $x$ to $y$. In particular,
when $P$ is connected, let $\sigma _{P}=\min_{e\in E\left( P\right) }\mathbf{%
\omega }\left( e\right) $, then 
\begin{equation*}
\left\vert Z\left( x\right) -Z\left( y\right) \right\vert \leq \frac{1}{%
\alpha }\left( \sigma _{P}\right) ^{\alpha -1}\left\Vert x-y\right\Vert
\end{equation*}%
for any $x,y$ on the support of $G_{P}$.
\end{proposition}

\subsection{Schiffman's unified approach of landscape function}
To unify these two approaches, in \cite{Schiffman2021}, Schiffman proposed the following idea: Consider the first-variation of ${\bf M}_\alpha$ between real rectifiable 1-currents with finite ${\bf M}_\alpha$ mass. 
Then, define landscape function on a larger family: $\alpha$-balanced real rectifiable 1-currents.\\

\begin{definition} \label{def:var-T character}
 Given two rectifiable 1-currents $T = \rect(M,\theta,\xi), S= \rect(N,\varphi,\zeta) $ with finite ${\bf M}_\alpha$ mass
for $\alpha \in [0,1)$, the  $T$-character of $S$  is defined by
\begin{equation} \label{eqn:var, character def}
	\chi_{_T}(S) := 
 \begin{cases}
     \int_N \theta^{\alpha-1}\varphi d\haus, &\text{ if } \haus(N \setminus M)=0;\\
      \infty, &\text{ otherwise. }
 \end{cases}
\end{equation}
By setting $\theta(x)=0$ for all $x\in M^c$, Equation (\ref{eqn:var, character def}) can be simply expressed as
\[\chi_{_T}(S)=\int_N \theta^{\alpha-1}\varphi d\haus.\]

\end{definition}

We now consider the first-variation of ${\bf M}_\alpha$, and get the following result:\\
 \begin{theorem}[B. Schiffman, Q. Xia]
Suppose $T = \rect(M,\theta,\xi)$ and $S = \rect(N,\varphi,\zeta)$ are 
$\malpha$-finite rectifiable 1-currents for $\alpha \in [0,1)$.  If $\chi_{_T}(S) < +\infty$, then
\begin{equation}\label{eqn:var-first variation}
\frac{d}{d\varepsilon} \big[ \malpha(T + \varepsilon S) \big]
			\Big\rvert_{\varepsilon=0}
	= \alpha \int_N \theta^{\alpha-1}\varphi \langle \xi,\zeta \rangle d\haus.
\end{equation}
\end{theorem}
\begin{proof} Since $\chi_{_T}(S)<\infty$, we have $\haus(N \setminus M) = 0$ and $\int_N \theta^{\alpha-1}\varphi d\haus<\infty$. Direct calculation gives
\[
	\frac{d}{d\varepsilon} \big[ \malpha(T + \varepsilon S) \big] 
			\Big\rvert_{\varepsilon=0} 
	=\lim\limits_{\varepsilon \rightarrow 0} \int_{N} \frac{|\theta(x)+\varepsilon\varphi(x)\langle \xi,\zeta \rangle|^\alpha-\theta(x)^\alpha } {\varepsilon}d\haus,
	\]
provided the limit exists.

By the strict concavity of $f(t)=t^\alpha$, for $\haus$-a.e. $x\in N$, 
 \[\left\lvert\frac{|\theta(x)+\varepsilon\varphi(x)\langle \xi,\zeta \rangle|^\alpha-\theta(x)^\alpha } {\varepsilon}\right\rvert\leq \alpha  \theta(x)^{\alpha-1} \varphi(x).\]
Thus,
	\[\int_{N} \left\lvert \frac{|\theta(x)+\varepsilon\varphi(x)\langle \xi,\zeta \rangle|^\alpha-\theta(x)^\alpha } {\varepsilon} \right\rvert d\haus
	\leq \int_{N} \theta(x)^{\alpha-1}\varphi(x) d\haus = \chi_{_T}(S)
	\]
for all nonzero $\varepsilon$.
Using dominated convergence to exchange the limit and integration, 
\begin{eqnarray*}
	\frac{d}{d\varepsilon} \big[ \malpha(T + \varepsilon S) \big] 
			\Big\rvert_{\varepsilon=0} 
	&=&  \int_{N}\lim\limits_{\varepsilon \rightarrow 0} \frac{|\theta+\varepsilon\varphi\langle \xi,\zeta \rangle|^\alpha-\theta^\alpha } {\varepsilon}d\haus\\
&	=& \alpha\int_N \theta^{\alpha-1}\varphi \langle \xi,\zeta \rangle d\haus.
\end{eqnarray*}
\end{proof}

Note that each oriented Lipschitz curve $\gamma$ associates with a rectifiable 1-current $I_\gamma$. For a given  $T = \rect(M,\theta,\xi)$, set
$
	\chi_{_T}(\gamma) := \chi_{_T} \big( I_\gamma\big)
$.
When $\gamma$ is simple and is contained $\haus$-a.e.\ in $M$, 
	\[
	\chi_{_T}(\gamma) 
	= \int_\gamma \theta^{\alpha-1}  d\haus
.	\]
\begin{definition} \label{def:balanced system} \label{def:var-balanced system}
Given $\alpha \in [0,1]$, a rectifiable 1-current $T = \rect(M,\theta,\xi)$ is called $\alpha$-balanced 
if
\begin{equation} \label{eqn:balanced defn}
	\int_C \theta^{\alpha-1}\langle \xi,\tau_{_C} \rangle d\haus = 0
\end{equation}
for every oriented Lipschitz closed curve $C$ satisfying 
$$\chiT(C) :=\int_C \theta^{\alpha-1} d\mathcal{H}^1 < +\infty,$$ 
where $\tau_{_C}$ is the unit-tangent orientation of $C$.
\end{definition}

We now define the  landscape functions of $\alpha$-balanced rectifiable 1-currents as follows:

\begin{definition} \label{def:var-landscape fns}
Let $T = \rect(M,\theta,\xi)$ be an $\alpha$-balanced rectifiable 1-current.
Define the $\alpha-$landscape function 
$Z_{_T} : \mathbb R^d \times \mathbb R^d \rightarrow (-\infty, + \infty]$ 
by
	\[
	Z_{_T}(p,q) 
	:=\int_\gamma \theta^{\alpha-1} \langle \xi,\tau_\gamma \rangle d\haus 
	\]
for any Lipschitz, $\chi_{_T}$-finite curve $\gamma$ from $p$ to $q$ and $Z_{_T}(p,q) 
	:=+ \infty $ if no such $\gamma$ exists.
\end{definition}

We now describe some examples of $\alpha$-balanced currents.

The first type of examples are $\alpha$-efficient rectifiable currents.
Recall that an $\alpha$-optimal  rectifiable 1-current $T$ 
is an $\mathbf{M}_{\alpha}$ minimizer among all 1-currents with the same boundary 
(those achieving the same transport): 
\[\mathbf M_\alpha(T) \leq \mathbf{M}_\alpha(T + \Ccal)\]
for every rectifiable 1-current $\Ccal$ with $\partial \Ccal=0$.

Now, we use the concept of $T$-character to introduce the following modification:
\begin{definition} \label{def:BandE-efficient systems}
A rectifiable 1-current $T $ is called $\alpha$-efficient if 
\[\mathbf M_\alpha(T) \leq \mathbf{M}_\alpha(T + \Ccal)\] for every 
rectifiable 1-current $\Ccal$ with $\partial \Ccal=0$ and $\chiT(\Ccal)<+\infty$.
\end{definition}
Immediately from this definition, there are two major families 
of $\alpha$-efficient rectifiable 1-currents:
\begin{itemize}
\item For a given $\alpha \in [0,1)$, every $\alpha$-optimal 
	rectifiable 1-current is $\alpha$-efficient.
\item Every finite acyclic weighted directed graph $G$ is $\alpha$-efficient 
	for all $\alpha \in [0,1)$ as there exist no nonzero cycles $\Ccal$ 
	that  are $\chiT[G]$-finite.
\end{itemize}

\begin{example}
\begin{figure}[h]
\includegraphics[width=2in]{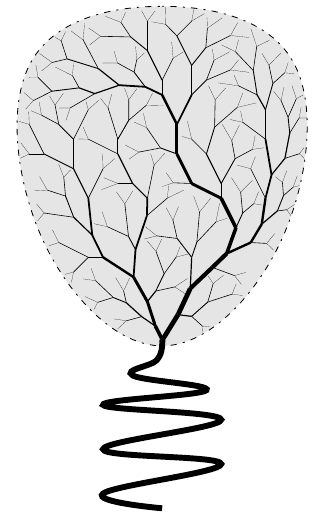} 
\caption{An example of a non-trivial $\alpha$-efficient rectifiable current}\label{egBalloon_fig}
\end{figure}
Figure \ref{egBalloon_fig} gives another example of $\alpha$-efficient rectifiable currents.
In this figure, the upper portion on its 
own is $\tilde \alpha$-optimal for a specific $\tilde \alpha \in (0,1)$, 
but the bottom portion is not optimal for any $\alpha$, thus the entire system together is 
not optimal for $\tilde \alpha$ or any other value in $[0,1]$.
However there are no finite-character cycles that lower the
$\malpha[\tilde \alpha]$-cost of 
the top or bottom portions, hence this system is $\tilde \alpha$-efficient.

\end{example}
\begin{theorem}\label{prop:BandE-efficient currents variations}
Let $T $ be $\alpha$-efficient with ${\bf M}_{\alpha}(T)<\infty$.
\begin{itemize}
	\item[{a.)}] \label{item:BandE-zero boundary variations}
	If $\Ccal$ has finite $T$-character and $\partial \Ccal = 0$ then
	\[
	\frac{d}{d\varepsilon} \big [ \mathbf{M}_\alpha(T+ \varepsilon \Ccal) 
		\big ]\Big\rvert_{\varepsilon=0} = 0.
	\]
\item[{b.)}] \label{item:BandE-equal boundary variations}
	If $S$ has finite $T$-character and $\partial S = \partial T$ then
	\[
	\frac{d}{d\varepsilon} \big [ \mathbf{M}_\alpha(T+ \varepsilon S) 
		\big ]\Big\rvert_{\varepsilon=0} = \alpha \mathbf{M}_\alpha(T).
	\]
\end{itemize}
\end{theorem}

\begin{corollary}
Suppose $T$ is $\alpha$-efficient for some $\alpha\in (0,1)$, then $T$ is $\alpha$-balanced.
\end{corollary}

Since every $\alpha$-efficient rectifiable 
1-current is also $\alpha$-balanced, one may consider the corresponding $\alpha$-landscape function using Definition
\ref{def:var-landscape fns}.
Observe that both definitions by Santambrogio and Xia  agrees with the one using Definition
\ref{def:var-landscape fns}.
Also noting some $\alpha$-efficient systems may be neither finite 
nor $\alpha$-optimal, consequently: the definition of landscape function for $\alpha$-balanced currents
 unifies and generalizes the existing definitions 
of landscape function provided by Santambrogio and Xia.

We now describe another type of examples for $\alpha$-balanced currents: Finite directed graphs containing cycles. 
\begin{example}
Given a finite directed graph $G = [V(G), E(G), w]$, let $Q_G$ be the corresponding vertex-edge incidence matrix, and 
\[\vec{M}_\alpha:=[w(e_1)^{\alpha-1}\mathcal{H}^1(e_1),\cdots, w(e_n)^{\alpha-1}\mathcal{H}^1(e_n)].\]
Then, using definition, $G$ is $\alpha$-balanced when
\[Q_G\vec{c}=\vec{0}\Rightarrow \vec{M}_\alpha\vec{c}=\vec{0}.,\quad i.e., ker(Q_G)\subseteq ker(\vec{M}_\alpha).\]

It turns out that:
A finite directed graph $G = [V(G), E(G), w]$ is
$\alpha$-balanced if and only if
\begin{equation} \label{eqn:BandE-finite cycles}
 \sum_{e \in E(\mathcal C)} w(e)^{\alpha-1}  sgn(c(e)) \haus(e)
	= 0
\end{equation}
for all simple cycles $\mathcal C= \sum_{e \in E(G)} c(e) [[e]]$ contained in $G$. 
Here, simple means $sgn(c(e)) \in \{-1, 0 , 1\}$ for all $e$.

\end{example}

\begin{theorem} [B. Schiffman, Q. Xia]
Let $T\in Path(\mu, \nu)$ be $\alpha$-balanced. If $T$ is compatible with a transport plan $\pi \in Plan(\mu,\nu)$, 
then
\begin{equation} \label{eqn:malpha-pi integral identity}
	\malpha(T) 
	= \iint_{\mathbb R^d \times \mathbb R^d} Z_T(x,y) d\pi.
\end{equation}
\end{theorem}	
\begin{proof}[Proof Sketch:]
\begin{itemize}
\item Express $Z_T(x,y)=z(y)-z(x)$ with $\malpha(T) =\int_{\mathbb R^d} z d(\partial T)$.
\item Compatibility implies that $Z_T(x,y)\ge 0$ for $\pi$-a.e. $(x,y)$.
\item Show that $\malpha(T) 
	= \iint_{\mathbb R^d \times \mathbb R^d} Z_T d\pi$.
\end{itemize}
\end{proof}

\begin{remark}
    Note that, (\ref{eqn:malpha-pi integral identity}) expresses the $\mathbf{M}_\alpha$-cost of $T$ in a format similar to the Monge-Kantorovich transport cost. This observation may deserve further exploration.
\end{remark}

\section{Open questions and future directions}

\textcolor{blue}{We conclude by mentioning a few directions in which the ideas described in this article may be further developed. These questions are not meant to be exhaustive, but rather to indicate some ways in which ramified optimal transportation continues to interact with geometric measure theory, numerical analysis, and the modeling of branching structures in nature.}

\textcolor{blue}{A first direction concerns the relation between ramified optimal transportation and Gilbert's original minimum-cost communication network problem \cite{Gilbert1967}. As explained above, the transport-path formulation studied in this article is a directed model with prescribed boundary $\partial G=b-a$. This makes it natural to regard the problem as a Plateau-type problem for weighted polyhedral $1$-chains, and later for rectifiable currents. Gilbert's original problem is different in structure. Even in the discrete setting, it involves pairwise channel demands between finitely many cities and undirected or bidirectional network links, rather than a directed weighted graph connecting prescribed source and target measures. A natural question is whether one can develop a geometric-measure-theoretic formulation closer to Gilbert's original model. In the discrete case, such a formulation would need to encode a matrix of pairwise demands, the sharing of common infrastructure, and the corresponding economy-of-scale cost. A further challenge is to pass from finitely many cities to general Radon measures. In such a continuum version, the communication demand would no longer be described by a finite matrix, but rather by a measure on pairs of points, or by a suitable traffic plan. The network itself would then need to carry shared bidirectional traffic while retaining enough geometric structure to allow variational analysis. It would be interesting to understand whether such a model admits existence theorems, regularity properties, approximation schemes, or current-theoretic or varifold-like descriptions of minimizers.}

\textcolor{blue}{A second direction is to develop more sophisticated biological models based on ramified transport. The leaf model described in Section~2 captures one essential mechanism: the competition between the benefit of generating new tissue and the transportation cost required to support it. However, real plant growth involves many additional effects, including anisotropic growth, mechanical constraints, chemical signaling, local environmental variation, and feedback between transport efficiency and tissue development. Similar questions arise for tree branches and root systems, where transport efficiency interacts with light exposure, nutrient acquisition, gravity, soil heterogeneity, and mechanical stability. A natural problem is to couple ramified transport with these biological and physical mechanisms while preserving a tractable variational structure. Such models could help explain not only the shape of individual leaves, but also the emergence of branching architectures in roots, stems, and whole plants. The same framework may also be useful for modeling degeneration. In the growth model, increasing values of the environmental parameter $\epsilon$ allow more cells to be generated. Conversely, one may consider a decreasing sequence of $\epsilon$ values, representing a reduction in environmental support or metabolic benefit. Cells whose transportation expense is no longer justified by the reduced revenue may then be removed from the leaf. In this way, the model could describe not only the formation of a leaf, but also its gradual decay. Since different regions of a leaf may cease to be supported at different stages, such a decreasing-revenue process may provide a mathematical mechanism for describing spatial patterns of leaf coloration during senescence.}

\textcolor{blue}{A third direction concerns numerical approximation. The multiscale subdivision--refinement algorithm described in Section~1.3 gives an effective way to produce high-quality approximations of ramified optimal transport paths for large atomic measures. Its practical success comes from combining exact optimization on small subproblems with energy-decreasing local and global refinements. Nevertheless, for large point sets the method remains heuristic: it does not provide a certificate of global optimality, and its convergence properties are not fully understood. It would be useful to develop numerical methods with stronger theoretical guarantees, such as convergence to minimizers under mesh refinement, computable lower bounds, a posteriori error estimates, or optimality certificates for restricted classes of networks. Another related question is how to compare different numerical formulations of branched transport, including Eulerian, Lagrangian, graph-based, and current-theoretic approaches. Progress in this direction would not only improve numerical simulations, but also clarify which geometric features of optimal transport paths are stable under approximation.}

\textcolor{blue}{These directions illustrate the broader theme of the article: branching patterns can be studied simultaneously as geometric objects, variational minimizers, transport networks, and models of natural growth. Developing these connections further may lead to new mathematical tools for understanding branching structures in both living and nonliving systems.}

\bibliographystyle{plain}
\bibliography{refs}
\end{document}